\documentclass[a4paper,12pt]{article}

\usepackage[applemac]{inputenc}
\usepackage{geometry}
\usepackage[final]{pdfpages}
\usepackage[english,french]{babel}
\usepackage{paralist}
\usepackage{pifont}
\usepackage{amsmath,amsfonts,amssymb,amsthm,cases}
\usepackage{empheq}
\usepackage{array,dcolumn,booktabs}
\usepackage[section]{placeins}
\usepackage{tikz,todonotes}
\usepackage{subfig}
\usepackage{url}
\usepackage{microtype}
\usepackage{marvosym}
\usepackage{graphicx}
\usepackage{hyperref}
\usepackage{epstopdf}
\usepackage{makeidx}
\usepackage{multicol}
\graphicspath{{images/}}
\usepackage{listings}
\def\twoplot[#1]#2#3#4#5{
\begin{figure}[hbt]
\begin{multicols}{2}
\begin{center}
    \includegraphics*[#1]{#2}
    \caption{\label{#2} #4}
\end{center}
\begin{center}
    \includegraphics*[#1]{#3}
    \caption{\label{#3} #5}
\end{center}
\end{multicols}
\end{figure}
}
\def\threeplot[#1][#2][#3]#4#5#6#7#8#9{
\begin{figure}[hbt]
\begin{multicols}{3}
\begin{center}
    \includegraphics*[#1]{#4}
    \caption{\label{#4} #7}
\end{center}
\begin{center}
    \includegraphics*[#2]{#5}
    \caption{\label{#5} #8}
\end{center}
\begin{center}
    \includegraphics*[#3]{#6}
    \caption{\label{#6} #9}
\end{center}
\end{multicols}
\end{figure}
}

\newcommand{\ds}{\displaystyle}
\newcommand{\bg}{\begin{equation}}
\newcommand{\ed}{\end{equation}}
\usepackage{xspace}
\makeatletter
\DeclareRobustCommand{\freefem}
{\valign{\vfil\hbox{##}\vfil\cr
   \textsf{FreeFem\kern-.1em}\cr
   $\hbox{\fontsize{\sf@size}{0}\textbf{+\kern-0.05em+}}$\cr}\xspace%
}
\hypersetup{
backref=true, 
pagebackref=true,
hyperindex=true, 
colorlinks=true, 
breaklinks=true, 
urlcolor= black, 
linkcolor= black, 
citecolor = black, 
 bookmarks=true, 
bookmarksopen=true, 
pdftitle={FreeFem++, a tool to solve PDEs numerically}, 
pdfauthor={Georges SADAKA}, 
pdfsubject={EDP 2D/3D FreeFem++} 
} 
\makeatother
\usepackage{textcomp}
\theoremstyle{plain}

\theoremstyle{definition}

\theoremstyle{remark}
\newtheorem*{rem}{Remark}

\def\p{\partial}
\def\n{\nabla}
\def\R{\mathbb{R}}

\begin{document}
{\selectlanguage{english} 
\title{ {\sc \textbf{\freefem, a tool to solve PDEs numerically}}}
\author{{\sc \Large{Georges Sadaka}}\\ \\
 \footnotesize{LAMFA CNRS UMR 7352}\\
\footnotesize{Université de Picardie Jules Verne}\\
\footnotesize{33, rue Saint-Leu, 80039 Amiens, France}\\
\small{\url{http://lamfa.u-picardie.fr/sadaka/}}\\
  \Email\xspace \small{georges.sadaka@u-picardie.fr}}

\maketitle

\begin{abstract}
\freefem is an open source platform to solve partial differential equations numerically, based on finite element methods. It was developed at the Laboratoire Jacques-Louis Lions, Université Pierre et Marie Curie, Paris by Frédéric Hecht in collaboration with Olivier Pironneau, Jacques Morice, Antoine Le Hyaric and Kohji Ohtsuka.\\

The \freefem platform has been developed to facilitate teaching and basic research through prototyping. \freefem has an advanced automatic mesh generator, capable of a posteriori mesh adaptation; it has a general purpose elliptic solver interfaced with fast algorithms such as the multi-frontal method UMFPACK, SuperLU . Hyperbolic and parabolic problems are solved by iterative algorithms prescribed by the user with the high level language of \freefem. It has several triangular finite elements, including discontinuous elements. For the moment this platform is restricted to the numerical simulations of problems which admit a variational formulation. \\

We will give in the sequel an introduction to \freefem which include the basic of this software. You may find more information throw this link \url{http://www.freefem.org/ff++}.
\end{abstract}

\tableofcontents
\clearpage
\section{Introduction}

\freefem is a \textbf{Free} software to solve PDE using the \textbf{F}inite \textbf{e}lement \textbf{m}ethod and it run on Mac, Unix and Window architecture.\\

In \freefem, it's used a user language to set and control the problem. This language allows for a quick specification of linear PDE's, with the variational formulation of a linear
steady state problem and the user can write they own script to solve non linear problem and time depend problem.\newline

It's a interesting tool for the problem of average size. It's also a help for the modeling in the sense where it allows to obtain quickly numerical results which is useful for modifying a physical model, to clear the avenues of Mathematical analysis investigation, etc ...\\

A documentation of \freefem is accessible on \url{www.freefem.org/ff++}, on the following link \url{www.freefem.org/ff++/ftp/FreeFem++doc.pdf}, you may also have a documentation in spanish on the following link \url{http://www.freefem.org/ff++/ftp/freefem++Spanish.pdf}\newline

You can also download an integrated environment called \freefem-cs, written by Antoine Le Hyaric on the following link \url{www.ann.jussieu.fr/~lehyaric/ffcs/install.php}
\section{Characteristics of FreeFem++}
Many of \freefem characteristics are cited in the full documentation of \freefem, we cite here some of them :
\begin{itemize}
\item Multi-variables, multi-equations, bi-dimensional and three-dimensional static or time dependent, linear or nonlinear coupled systems; however the user is required to describe the iterative procedures which reduce the problem to a set of linear problems.
\item Easy geometric input by analytic description of boundaries by pieces, with specification by the user of the intersection of boundaries.
\item Automatic mesh generator, based on the Delaunay-Voronoi algorithm \cite{LucPir98}\label{LucPir981}.
\item load and save Mesh, solution.
\item Problem description (real or complex valued) by their variational formulations, the write of the variational formulation is too close for that written on a paper.
\item Metric-based anisotropic mesh adaptation.
\item A large variety of triangular finite elements : linear, quadratic Lagrangian elements
and more, discontinuous P1 and Raviart-Thomas elements, ...
\item Automatic Building of Mass/Rigid Matrices and second member.
\item Automatic interpolation of data from a mesh to an other one, so a finite element function is view as a function of (x; y) or as an array.
\item LU, Cholesky, Crout, CG, GMRES, UMFPack sparse linear solver.
\item Tools to define discontinuous Galerkin finite element formulations P0, P1dc, P2dc
and keywords: {\ttfamily{\textcolor{red}{jump}, \textcolor{red}{mean}, \textcolor{red}{intalledges}}}.
\item Wide range of examples : Navier-Stokes, elasticity, fluid structure, eigenvalue
problem, Schwarz' domain decomposition algorithm, residual error indicator, ...
\item Link with other software : modulef, emc2, medit, gnuplot, ...
\item Generates Graphic/Text/File outputs.
\item A parallel version using mpi.
\end{itemize}
\section{How to start?}
All this information here are detailed in the \freefem documentation.
\subsection{Install}
First open the following web page $$\mbox{\url{http://www.freefem.org/ff++/}}$$

Choose your platform: Linux, Windows, MacOS X, or go to the end of the page to get the full list of downloads and then install by double click on the appropriate file.
\subsection{Text editor}
\begin{enumerate}
\item For Windows :\\
Install \textbf{notepad++} which is available at \url{http://notepad-plus.sourceforge.net/uk/site.htm}

\begin{itemize}

\item  Open Notepad++ and Enter  F5
\item   In the new window  enter  the command
\verb!FreeFem++ "$(FULL_CURRENT_PATH)"!
\item Click on  Save,  and enter  \texttt{FreeFem++} in the box  "Name", now  choose
the  short cut key  to launch directly FreeFem++  (for example \texttt{alt+shift+R})
 \item To add    Color Syntax Compatible with FreeFem++ In Notepad++,
 \begin{itemize}
\item In  Menu \verb!"Parameters"->"Configuration of the  Color Syntax"! proceed as follows:
\item In the list \verb!"Language"! select C++
\item Add "edp" in the  field  \verb!"add ext"!
\item Select  \verb!"INSTRUCTION WORD"! in the  list \verb!"Description"! and in the field
\verb!"supple! \verb!mentary key word"!, cut and past the following list:

\textcolor{blue}{P0 P1 P2 P3 P4 P1dc P2dc P3dc P4dc RT0 RT1 RT2 RT3 RT4 macro plot int1d int2d
solve movemesh adaptmesh trunc checkmovemesh on func buildmesh square Eigenvalue min max
imag exec LinearCG NLCG Newton BFGS LinearGMRES
catch try intalledges jump average mean load savemesh convect abs
sin cos tan atan asin acos cotan sinh cosh tanh cotanh atanh asinh acosh pow
exp log log10 sqrt dx dy endl cout}

\item Select "TYPE WORD" in the list "Description" and ... " "supplementary key word",
cut and past the following list

\textcolor{blue}{mesh real fespace varf matrix problem string border complex ifstream ofstream}

\item Click on \texttt{Save \& Close}.  Now nodepad++ is configured.
\end{itemize}
\end{itemize}
\item For MacOS :\\
Install \textbf{Smultron} which is available at \url{http://smultron.sourceforge.net}. It comes ready with color syntax for .edp file.  To teach it to launch \freefem files, do a "command B" (i.e. the menu Tools/Handle Command/new command) and create a command which does
 \begin{verbatim}
 /usr/local/bin/FreeFem++-CoCoa %%p
 \end{verbatim}

\item For Linux :\\
Install \textbf{Kate} which is available at  \url{ftp://ftp.kde.org/pub/kde/stable/3.5.10/src/kdebase-3.5.10.tar.bz2}\\

To personalize with color syntax for .edp file, it suffices to take those given by \textbf{Kate} for c++ and to add the keywords of \freefem. Then, download edp.xml and save it in the directory ".kde/share/apps/katepart/syntax".\\

We may find other description for other text editor in the full documentation of \freefem.

\end{enumerate}
\subsection{Save and run}
All \freefem code must be saved with file extension .edp and to run them you may double click on the file on MacOS or Windows otherwise we note that this can also be done in terminal mode by :
\texttt{FreeFem++ mycode.edp}
\section{Syntax and some operators}
\subsection{Data types}
In essence \freefem is a \index{compiler} compiler:  its language is typed, polymorphic, with exception and reentrant. Every variable must be declared of a certain type,  in a  declarative statement;
each statement are separated from the next by a semicolon ``\texttt{;}''.\\
Another trick is to \emph{comment in and out} by using the ``{\ttfamily{\textcolor{red}{//}}}'' as in C++. We note that, we can also comment a paragraph by using ``{\ttfamily{\textcolor{red}{/* paragraph */}}}'' and in order to make a break during the computation, we can use ``{\ttfamily exit(0);}''. \\
The variable {\ttfamily verbosity} changes the level of internal printing (0, nothing (unless there are syntax errors), 1 few, 10 lots, etc. ...), the default value is 2 and the variable {\ttfamily clock()} gives the computer clock.\\
The language allows the manipulation of basic types :\\

\begin{itemize}
\item current coordinates : {\ttfamily{\textcolor{red}{x}}}, {\ttfamily{\textcolor{red}{y}}} and {\ttfamily{\textcolor{red}{z}};}
\item current differentials operators : {\ttfamily{\textcolor{red}{dx}}}$=\ds\frac{\partial}{\partial x}$, {\ttfamily{\textcolor{red}{dy}}}$=\ds\frac{\partial}{\partial y}$, {\ttfamily{\textcolor{red}{dz}}}$=\ds\frac{\partial}{\partial z}$, {\ttfamily{\textcolor{red}{dxy}}}$=\ds\frac{\partial}{\partial xy}$, {\ttfamily{\textcolor{red}{dxz}}}$=\ds\frac{\partial}{\partial xz}$, {\ttfamily{\textcolor{red}{dyz}}}$=\ds\frac{\partial}{\partial yz}$, {\ttfamily{\textcolor{red}{dxx}}}$=\ds\frac{\partial}{\partial xx}$, {\ttfamily{\textcolor{red}{dyy}}}$=\ds\frac{\partial}{\partial yy}$ and {\ttfamily{\textcolor{red}{dzz}}}$=\ds\frac{\partial}{\partial zz}$;
\item integers, example : {\ttfamily{\textcolor{red}{int} a=1;}} 
\item reals, example : {\ttfamily{\textcolor{red}{real} b=1.;}} (don't forget to put a point after the integer number)
\item complex, example : {\ttfamily{\textcolor{red}{complex} c=1.+3\textcolor{red}{i};}} 
\item strings, example : {\ttfamily{\textcolor{red}{string} test="toto";}}
\item arrays with real component, example: {\ttfamily{\textcolor{red}{real}[\textcolor{red}{int}] V(n);}} where {\ttfamily{n}} is the size of {\ttfamily{V}},
\item arrays with complex component, example: {\ttfamily{\textcolor{red}{complex}[\textcolor{red}{int}] V(n);}}
\item matrix with real component, example: {\ttfamily{\textcolor{red}{real}[\textcolor{red}{int},\textcolor{red}{int}] A(m,n);}}
\item matrix with complex component, example: {\ttfamily{\textcolor{red}{complex}[\textcolor{red}{int},\textcolor{red}{int}] C(m,n);}}
\item bidimensional (2D) finite element meshes, example : {\ttfamily{\textcolor{red}{mesh} Th;}}
\item 2D finite element spaces, example : {\ttfamily{\textcolor{red}{fespace} Vh(Th,P1);} \textcolor{red}{// where Vh is the Id space}}
\item threedimensional (3D) finite element meshes, example : {\ttfamily{\textcolor{red}{mesh3} Th3;}}
\item 3D finite element spaces, example : {\ttfamily{\textcolor{red}{fespace} Vh3(Th3,P13d);}} 
\item {\ttfamily{\textcolor{red}{int1d}(Th,$\Gamma$)( u*v )}} $= \displaystyle\int_\Gamma u\cdot v\,dx$ where $\Gamma\subset \mathbb{R}$;
\item {\ttfamily{\textcolor{red}{int2d}(Th)( u*v )}} $= \displaystyle\int_\Omega u\cdot v\,dxdy$ where $\Omega\subset \mathbb{R}^2$;
\item {\ttfamily{\textcolor{red}{int3d}(Th)( u*v )}} $= \displaystyle\int_\Omega u\cdot v\,dxdydz$ where $\Omega\subset \mathbb{R}^3$.

\end{itemize}
\subsection{Some operators}
We cite here some of the operator that are defined in \freefem :
\begin{lstlisting}[firstnumber=last]
+, -, *, /, ^, 
<, >, <=, >=, 
&, |, // where a & b = a and b, a | b = a or b
=, +=, -=, /=, *=, !=, ==.
\end{lstlisting}
\subsection{Manipulation of functions}
We can define a function as :
\begin{itemize}
\item an analytical function, example : {\ttfamily{\textcolor{red}{func} u0=exp(-x\small{$^\wedge$}2-y\small{$^\wedge$}2)},u1=1.*(x>=-2 \& x<=2)}; 
\item a finite element function or array, example : {\ttfamily{Vh u0=exp(-x\small{$^\wedge$}2-y\small{$^\wedge$}2)};}.\\
We note that, in this case {\ttfamily{u0}} is a finite element, thus {\ttfamily{u0[]}} return the values of {\ttfamily{u0}} at each degree of freedom and to have access to the $i^{{\small\mbox{th}}}$ element of {\ttfamily{u0[]}} we may use {\ttfamily{u0[][i]}}.\\
We can also have an access to the value of {\ttfamily{u0}} at the point {\ttfamily{(a,b)}} by using {\ttfamily{u0(a,b)}};
\item a complex value of finite element function or array, example : {\ttfamily{Vh<\textcolor{red}{complex}> u0=x+1\textcolor{red}{i}*y}};
\item a formal line function, example : {\ttfamily{\textcolor{red}{func real} g(\textcolor{red}{int} a, \textcolor{red}{real} b) \{ .....; \textcolor{red}{return} a+b;\}}} and to call this function for example we can use {\ttfamily{g(1,2)}}.\\
We can also put an array inside this function as :
\begin{lstlisting}[firstnumber=last]
func real f(int a, real[int] U){
	Vh NU;
	NU[]=U;
	return a*NU;
}
Vh U=x,FNU=f(5,U[]);
\end{lstlisting}
\item {\ttfamily{\textcolor{red}{macro}}} function, example : {\ttfamily{\textcolor{red}{macro} F(t,u,v)[t*\textcolor{red}{dx}(u),t*\textcolor{red}{dy}(v)]\textcolor{red}{//}}}, notice that every {\ttfamily{\textcolor{red}{macro}}} must end by ``{\ttfamily{\textcolor{red}{//}}}'', it simply replaces the name  {\ttfamily{F(t,u,v)}} by {\ttfamily{[t*\textcolor{red}{dx}(u),t*\textcolor{red}{dy}(v)]}} and to have access only to the first element of {\ttfamily{F}}, we can use  {\ttfamily{F(t,u,v)[0]}}. 

In fact, we note that the best way to define a function is to use {\ttfamily{\textcolor{red}{macro}}} function since in this example {\ttfamily{t,u}} and {\ttfamily{v}} could be integer, real, complex, array or finite element, ... \\
For example, here is the most used function defined by a {\ttfamily{\textcolor{red}{macro}}} :
\begin{lstlisting}[firstnumber=last]
macro Grad(u)[dx(u),dy(u)]// in 2D
macro Grad(u)[dx(u),dy(u),dz(u)]// in 3D
macro div(u,v)[dx(u)+dy(v)]// in 2D
macro div(u,v,w)[dx(u)+dy(v)+dz(w)]// in 3D
\end{lstlisting}

\end{itemize}

\subsection{Manipulation of arrays and matrices}
Like in matlab, we can define an array such as : {\ttfamily{\textcolor{red}{real}[\textcolor{red}{int}] U=1:2:10;}} which is an array of 5 values {\ttfamily{U[i]=2*i+1; i=0 to 4}} and to have access to the $i^{{\small\mbox{th}}}$ element of {\ttfamily{U}} we may use {\ttfamily{U(i)}}.\\

Also we can define a matrix such as {\ttfamily{\textcolor{red}{real}[\textcolor{red}{int},\textcolor{red}{int}] A=[ [1,2,3] , [2,3,4] ];}} which is a matrix of size $2\times 3$ and to have access to the $(i,j)^{{\small\mbox{th}}}$ element of {\ttfamily{A}} we may use {\ttfamily{A(i,j)}}.\\

We will give here some of manipulation of array and matrix that we can do with \freefem :
\begin{lstlisting}[firstnumber=last]
real[int] u1=[1,2,3],u2=2:4; // defining u1 and u2
real u1pu2=u1'*u2; // give the scalar product of u1 and u2, here u1' is the transpose of u1;
real[int] u1du2=u1./u2; // divided term by term
real[int] u1mu2=u1.*u2; // multiplied term by term
matrix A=u1*u2'; // product of u1 and the transpose of u2
matrix<complex> C=[ [1,1i],[1+2i,.5*1i] ];
real trA=trace([1,2,3]*[2,3,4]'); // trace of the matrix
real detA=det([ [1,2],[-2,1] ]); // just for matrix 1x1 and 2x2
\end{lstlisting}
\subsection{Loops and conditions}
The {\ttfamily{\textcolor{red}{for}}} and {\ttfamily{\textcolor{red}{while}}} loops are implemented in \freefem together with {\ttfamily{\textcolor{red}{break}}} and {\ttfamily{\textcolor{red}{continue}}} keywords.\\
In {\ttfamily{\textcolor{red}{for}}}-loop, there are three parameters; the INITIALIZATION of a control variable, the CONDITION to continue, the CHANGE of the control variable. While CONDITION is true, {\ttfamily{\textcolor{red}{for}}}-loop continue.
\begin{lstlisting}[firstnumber=last]
for (INITIALIZATION; CONDITION; CHANGE)
     { BLOCK of calculations }
\end{lstlisting}
An example below shows a sum from 1 to 10 with result is in {\ttfamily sum},
\begin{lstlisting}[firstnumber=last]
int sum=0;
for (int i=1; i<=10; i++)
   sum += i;
\end{lstlisting}
The while-loop
\begin{lstlisting}[firstnumber=last]
while (CONDITION) {
   BLOCK of calculations or change of control variables
}
\end{lstlisting}
is executed repeatedly until CONDITION become false.
The sum from 1 to 5 can also be computed by {\ttfamily{\textcolor{red}{while}}}, in this example, we want to show how we can exit from a loop in midstream by {\ttfamily \textcolor{red}{break}} and how the {\ttfamily{\textcolor{red}{continue}}} statement will pass the part from {\ttfamily{\textcolor{red}{continue}}} to the end of the loop :
\begin{lstlisting}[firstnumber=last]
int i=1, sum=0;
while (i<=10) {
   sum += i; i++;
   if (sum>0) continue;
   if (i==5) break;
}
\end{lstlisting}

\subsection{Input and output data}
The syntax of input/output statements is similar to C++ syntax. It uses {\ttfamily{\textcolor{red}{cout}}}, {\ttfamily{\textcolor{red}{cin}}}, {\ttfamily{\textcolor{red}{endl}}},
{\ttfamily{<<}} and {\ttfamily{>>}} :
\begin{lstlisting}[firstnumber=last]
int i;
cout << " std-out" << endl;
cout << " enter i= ? ";
cin >> i ;
Vh uh=x+y;
ofstream f("toto.txt");	f << uh[]; // to save the solution
ifstream f("toto.txt");	f >> uh[]; // to read the solution
\end{lstlisting}

We will present in the sequel, some useful script to use the \freefem data with other software such as \textbf{ffglut}, \textbf{Gnuplot}\footnote{\url{http://www.gnuplot.info/}}, \textbf{Medit}\footnote{\url{http://www.ann.jussieu.fr/~frey/software.html}}, \textbf{Matlab}\footnote{\url{http://www.mathworks.fr/products/matlab/}}, \textbf{Mathematica}\footnote{\url{http://www.wolfram.com/mathematica/}}, \textbf{Visit}\footnote{\url{https://wci.llnl.gov/codes/visit/}} when we save data with extension as \textbf{.eps}, \textbf{.gnu}, \textbf{.gp}, \textbf{.mesh}, \textbf{.sol}, \textbf{.bb},  \textbf{.txt} and \textbf{.vtu}. \\
For \textbf{ffglut} which is the visualization tools through a pipe of \freefem, we can plot the solution and save it with a .eps format such as :
\begin{lstlisting}[firstnumber=last]
plot(uh,cmm="t="+t+"  ;||u||_L^2="+NORML2[kk], fill=true,value=true,dim=2);
\end{lstlisting}
\begin{figure}[!htb]
\begin{center}

\includegraphics[width=13cm,height=9cm]{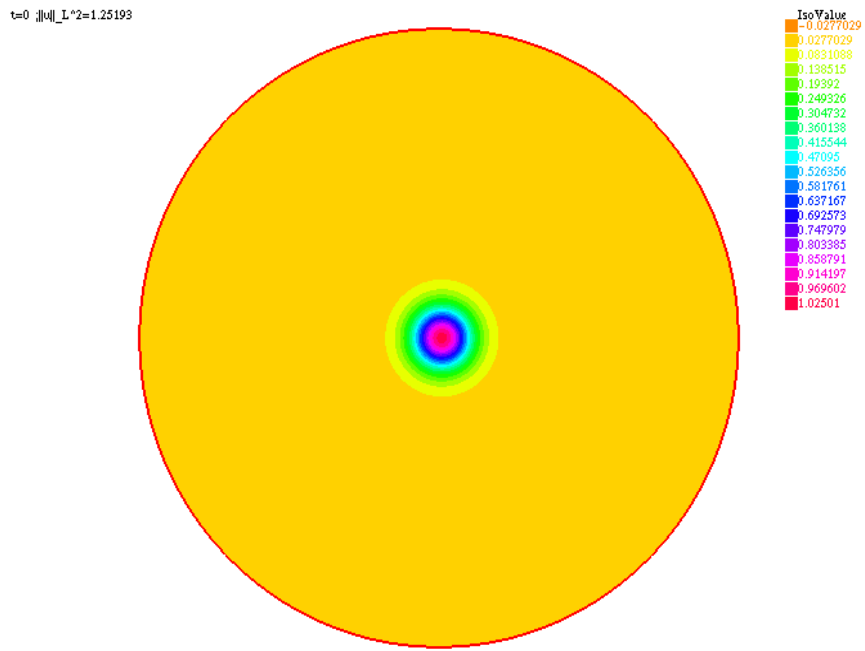}

\caption{\label{ffglut} Visualising of the solution using ffglut}
\end{center}
\end{figure}
For \textbf{Gnuplot}, we can save the data with extension \textbf{.gnu} or \textbf{.gp} such as :
\begin{lstlisting}[firstnumber=last]
{ ofstream gnu("plot.gnu"); // or plot.gp
//ofstream gnu("plot."+1000+k".gnu");  // to save the data
for (int i=0;i<=n;i++)
   gnu<<xx[i]<<" "<<yy[i]<<endl; // to plot yy[i] vs xx[i]
}
exec("echo 'plot \"plot.gnu\" w lp \
pause 5 \
quit' | gnuplot");
\end{lstlisting}

For \textbf{Medit}, we can save the data with extension \textbf{.mesh} and \textbf{.sol} such as :
\begin{lstlisting}[firstnumber=last]
load "medit"
int k=0;
savemesh(Th,"solution."+(1000+k)+".mesh");
savesol("solution."+(1000+k)+".sol",Th,uh);
medit("solution",Th,uh); // to plot the solution here
k+=1;
\end{lstlisting}
And then throw a terminal, in order to visualize the movie of the first 11 saved data, we can type this line :
\begin{lstlisting}[firstnumber=last]
ffmedit -a 1000 1010 solution.1000.sol
\end{lstlisting}
Don't forget in the window of \textbf{Medit}, to click on ``m'' to visualize the solution!\\
\clearpage
\begin{figure}[!htb]

\includegraphics[width=16cm,height=9cm]{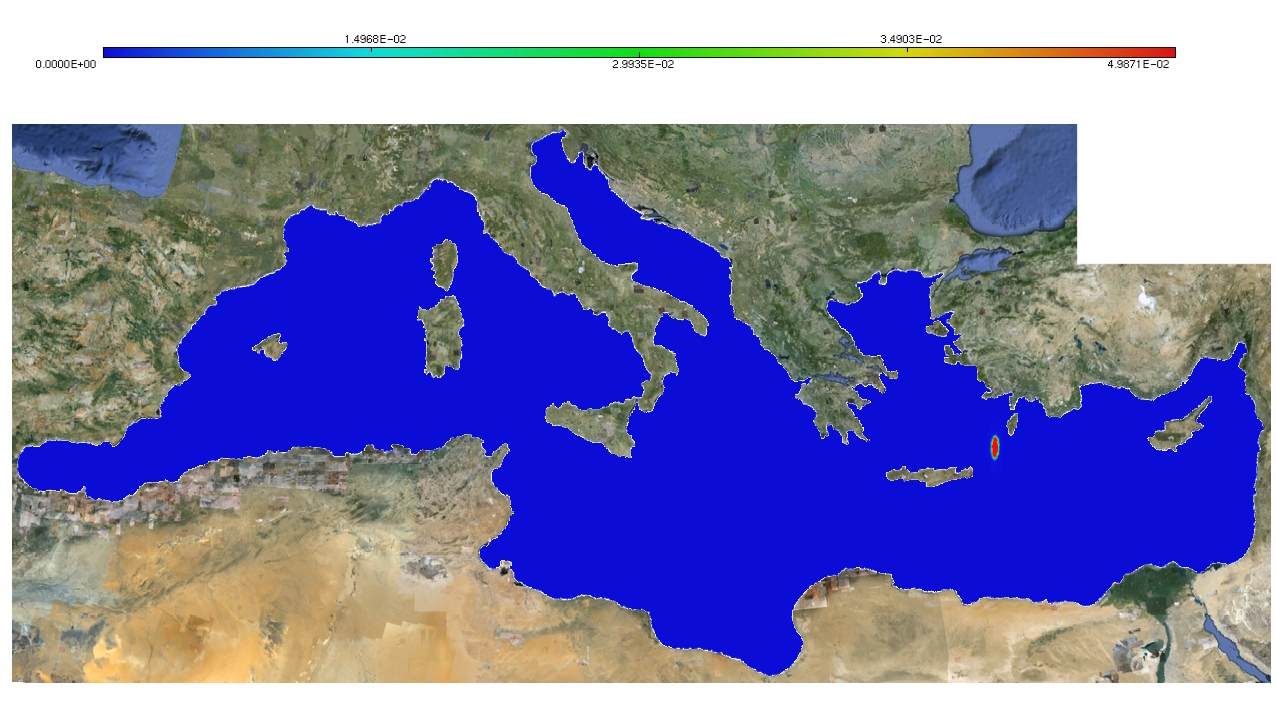}

\caption{\label{Medit} Visualising of the solution using Medit}
\end{figure}

For \textbf{Matlab}, we can save the data with extension \textbf{.bb} such as :
\begin{lstlisting}[firstnumber=last]
{ ofstream file("solution.bb"); 
   file << "2 1 1 "<< Vh.ndof << " 2 \n"; 
   for (int j=0;j<Vh.ndof ; j++)
      file << uh[][j] << endl;
}
\end{lstlisting}
And in order to visualize with \textbf{Matlab}, you can see the script made by Julien Dambrine at \url{http://www.downloadplex.com/Publishers/Julien-Dambrine/Page-1-0-0-0-0.html}.\\

\begin{figure}[!htb]
\begin{tabular}{>{\centering}m{8cm}>{\centering}m{8cm}}
 \includegraphics[width=8cm,height=6cm]{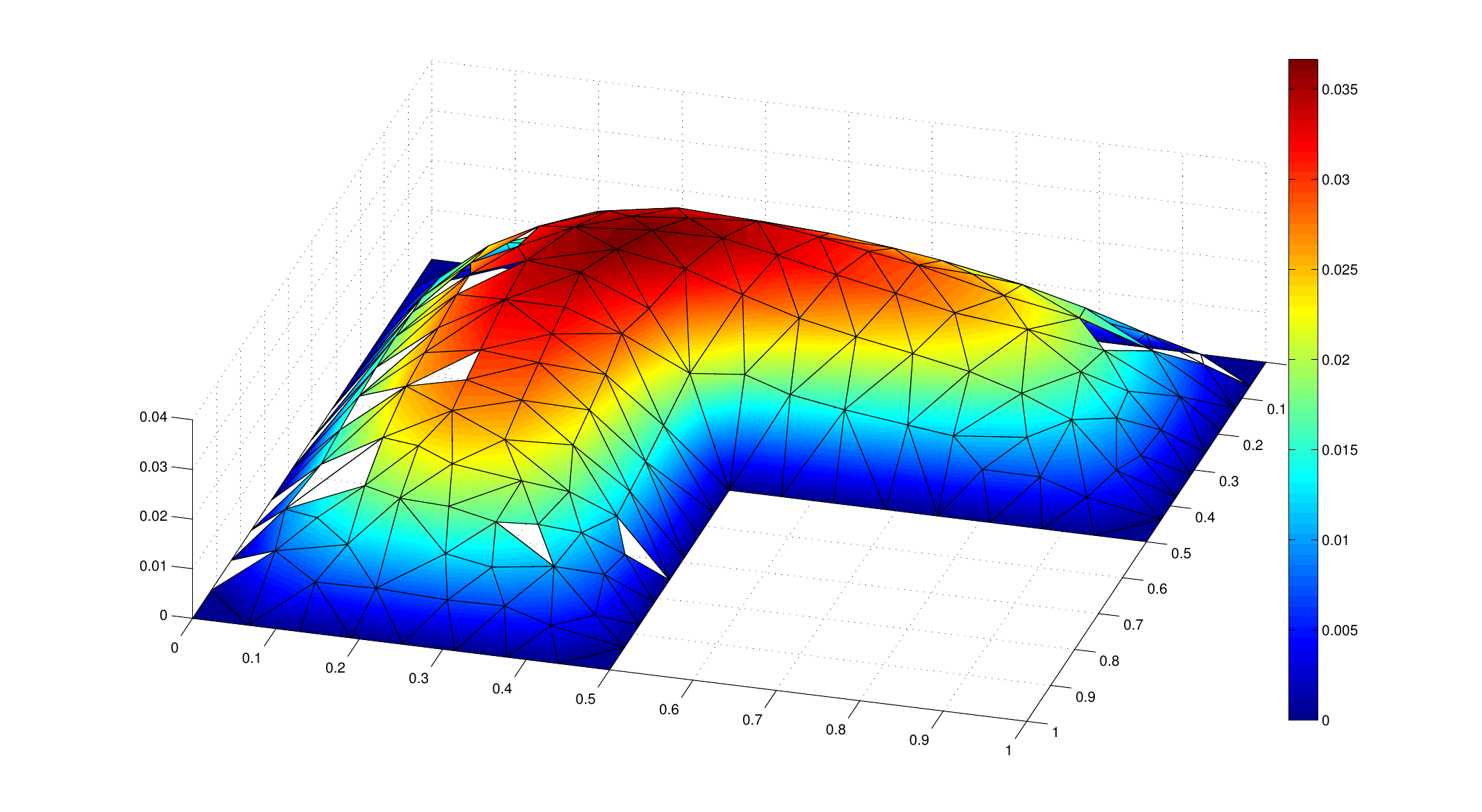}\includegraphics[height=6cm,width=8cm]{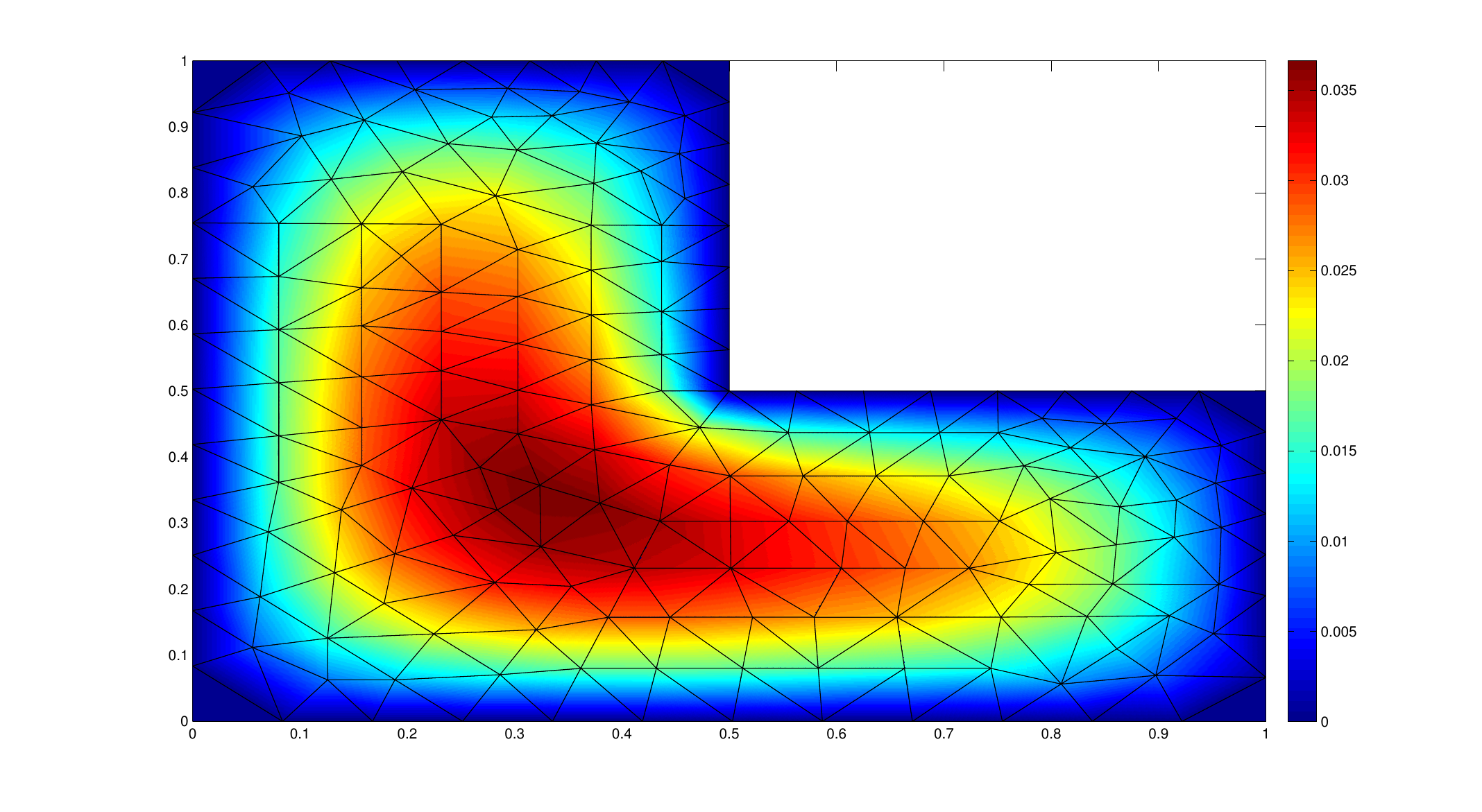}
\end{tabular}
\caption{\label{matlab} Visualising of the solution and the mesh using matlab}
\end{figure}
\clearpage
For \textbf{Mathematica}, we can save the data with extension \textbf{.txt} such as :
\begin{lstlisting}[firstnumber=last]
int k=0;
{ ofstream ff("uhsol."+(1000+k)+".txt");
for (int i=0;i<Th.nt;i++){
   for (int j=0; j <3; j++)
      ff<<Th[i][j].x<<" "<< Th[i][j].y<<" "<<uh[][Vh(i,j)]<<endl;
   ff<<Th[i][0].x<<" "<< Th[i][0].y<<" "<<uh[][Vh(i,0)]<<"\n";
}
}
k+=1;
\end{lstlisting}
\begin{figure}[!htb]
\begin{center}
 \includegraphics[width=14cm,height=8cm]{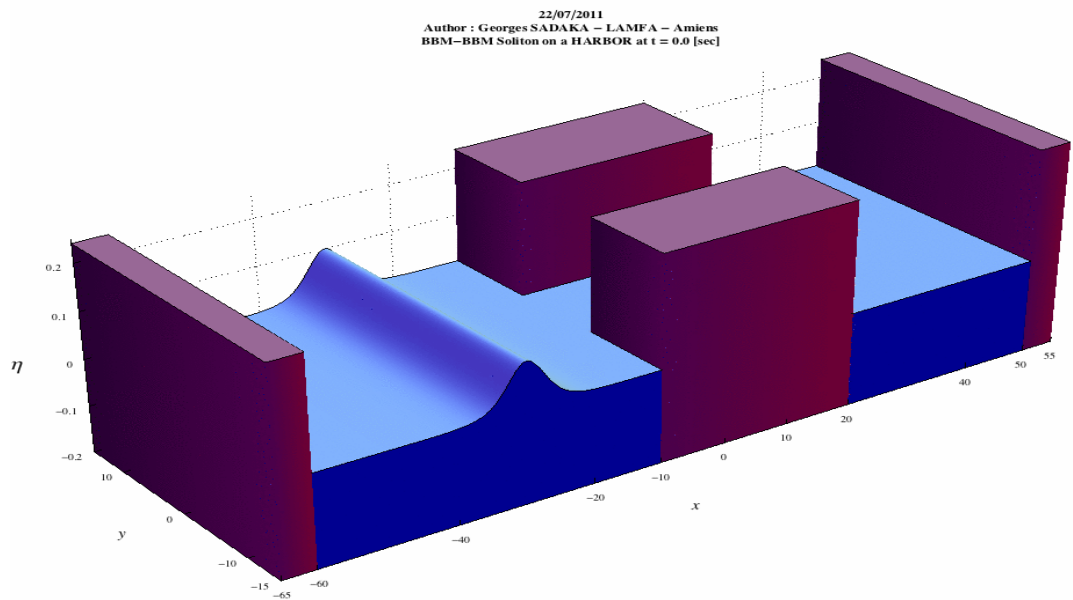}

\caption{\label{Mathematica} Visualising of the solution using Mathematica}
\end{center}
\end{figure}

For \textbf{Visit}, we can save the data with extension \textbf{.vtu} such as :
\begin{lstlisting}[firstnumber=last]
load "iovtk"
int k=0;
int[int] fforder2=[1,1,1];
savevtk("solution."+(1000+k)+".vtu",Th,uh1,uh2,order=fforder2,dataname="UH1 UH2", bin=true);
k+=1;
\end{lstlisting}
\begin{figure}[!htb]

 \includegraphics[width=16cm,height=8cm]{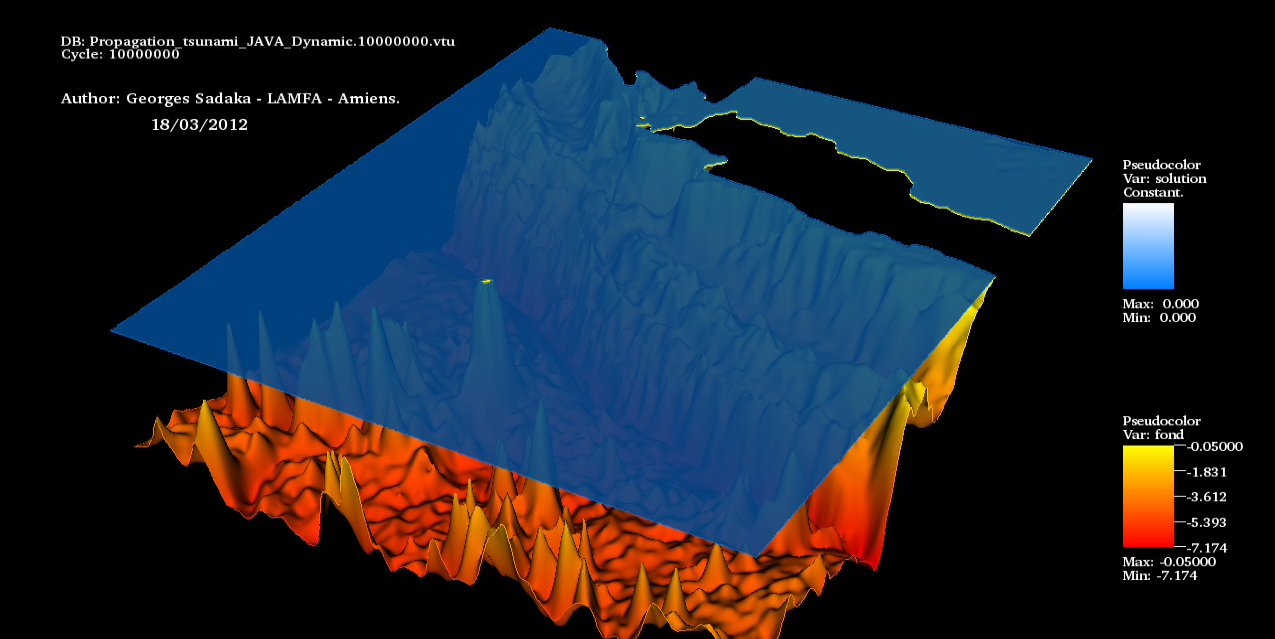}

\caption{\label{visit} Visualising of the solution using visit}
\end{figure}

\section{Construction of the domain $\Omega$}\label{consomega}
We note that in \freefem the domain is assumed to described by its boundary that is on the left side of the boundary which is implicitly oriented by the parametrization.\newline
Let $\Omega$ be the rectangle defined by its frontier $\partial \Omega=[-5,5]\times[-1,1]$ where his vertices are $A(-5,-1),B(5,-1),C(5,1)$ and $D(-5,1)$, so we must define the border $AB,BC,CD$ and $DA$ of $\partial\Omega$ by using the keyword \textcolor{red}{\ttfamily{border}} then the triangulation $\mathcal{T}_h$ of $\Omega$ is automatically generated by using the keyword \textcolor{red}{\ttfamily{buildmesh}}.

\begin{lstlisting}[firstnumber=last]
real Dx=.2; // discretization space parameter
int aa=-5,bb=5,cc=-1,dd=1;
border AB (t = aa, bb){x = t ;y = cc;label = 1;};
border BC (t = cc, dd){x = bb;y = t ;label = 2;};
border CD (t = bb, aa){x = t ;y = dd;label = 3;};
border DA (t = dd, cc){x = aa;y = t ;label = 4;};
mesh Th = buildmesh( AB(floor(abs(bb-aa)/Dx)) + BC(floor(abs(dd-cc)/Dx)) + CD(floor(abs(bb-aa)/Dx)) + DA(floor(abs(dd-cc)/Dx)) );
plot( AB(floor(abs(bb-aa)/Dx)) + BC(floor(abs(dd-cc)/Dx)) + CD(floor(abs(bb-aa)/Dx)) + DA(floor(abs(dd-cc)/Dx)) ); // to see the border
plot ( Th, ps="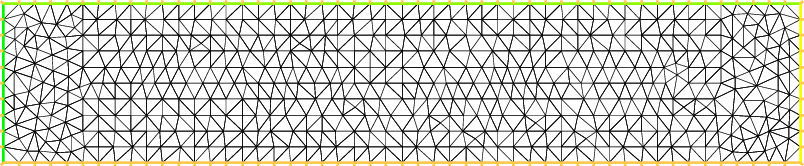"); // to see and save the mesh
\end{lstlisting}
The keyword \textcolor{red}{\ttfamily{label}} can be added to define a group of boundaries for later use (Boundary Conditions for instance). Boundaries can be referred to either by name ( $AB$ for example) or by label ( $1$ here).
\clearpage
\begin{figure}[!htb]
\begin{tabular}{>{\centering}m{8cm}>{\centering}m{8cm}}
\includegraphics[width=8cm,height=3cm]{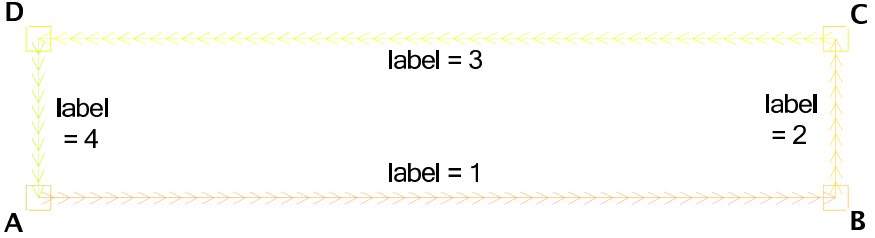} &
\includegraphics[width=8cm,height=2.3cm]{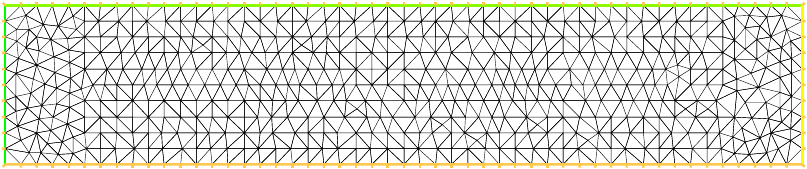}
\end{tabular}
\caption{\label{tab2} Plot of the border (left) and the mesh (right)}
\end{figure}
Another way to construct a rectangle domain with isotropic triangle is to use :
\begin{lstlisting}[firstnumber=last]
mesh Th=square(m,n,[x,y]); // build a square with m point on x direction and n point on y direction
mesh Th1=movemesh(Th,[x+1,y*2]); // translate the square ]0,1[*]0,1[ to a rectangle ]1,2[*]0,2[
savemesh(Th1,"Name.msh"); // to save the mesh
mesh Th2("mesh.msh"); // to load the mesh
\end{lstlisting}
\begin{figure}[htbp]
\begin{center}
  \includegraphics[height=6cm]{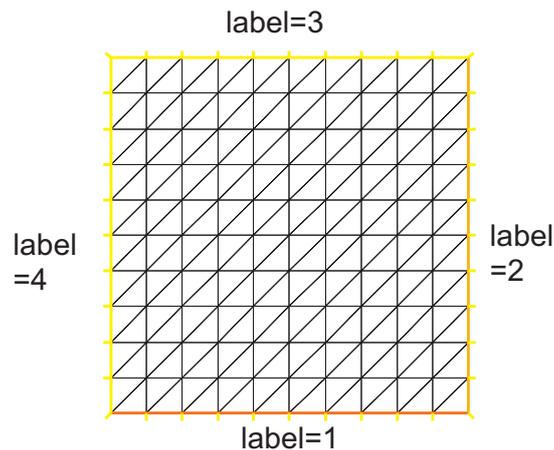}
\end{center}
  \caption{Boundary labels of the mesh by {\ttfamily{square(10,10)}}
  \label{fig:square}} \index{label}
\end{figure}
We can also construct our domain defined by a parametric coordinate as:
\begin{lstlisting}[firstnumber=last]
border C(t=0,2*pi){ x=cos(t);y=sin(t);label=1}
mesh Mesh_Name=buildmesh(C(50));
\end{lstlisting}
\clearpage
\begin{figure}[!htb]
\begin{center}
\includegraphics[height=5cm]{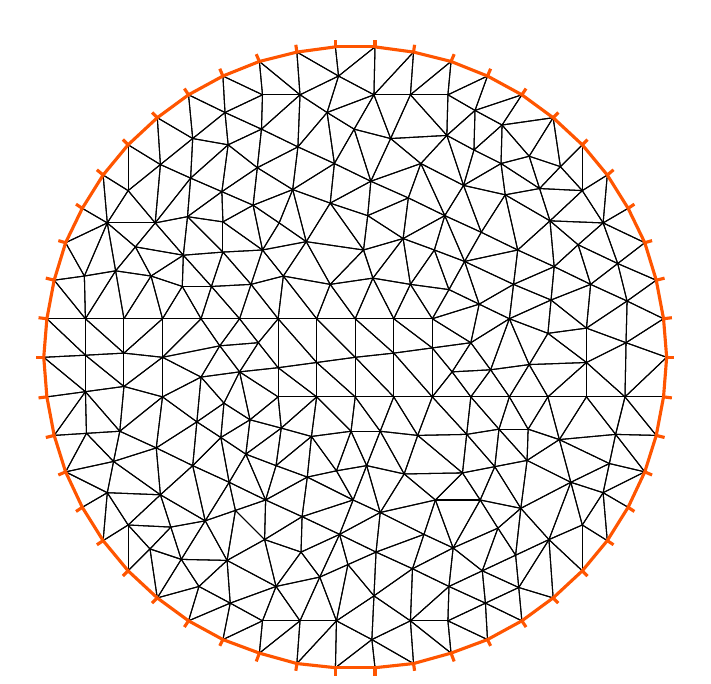}
\caption{\label{figmeshcirc} mesh {\ttfamily{Th}} by {\ttfamily{build(C(50))}}}
\end{center}
\end{figure}

To create a domain with a hole we can proceed as:
\begin{lstlisting}[firstnumber=last]
border a(t=0,2*pi){ x=cos(t); y=sin(t);label=1;}
border b(t=0,2*pi){ x=0.3+0.3*cos(t); y=0.3*sin(t);label=2;}
mesh Thwithouthole= buildmesh(a(50)+b(+30));
mesh Thwithhole   = buildmesh(a(50)+b(-30));
plot(Thwithouthole,wait=1,ps="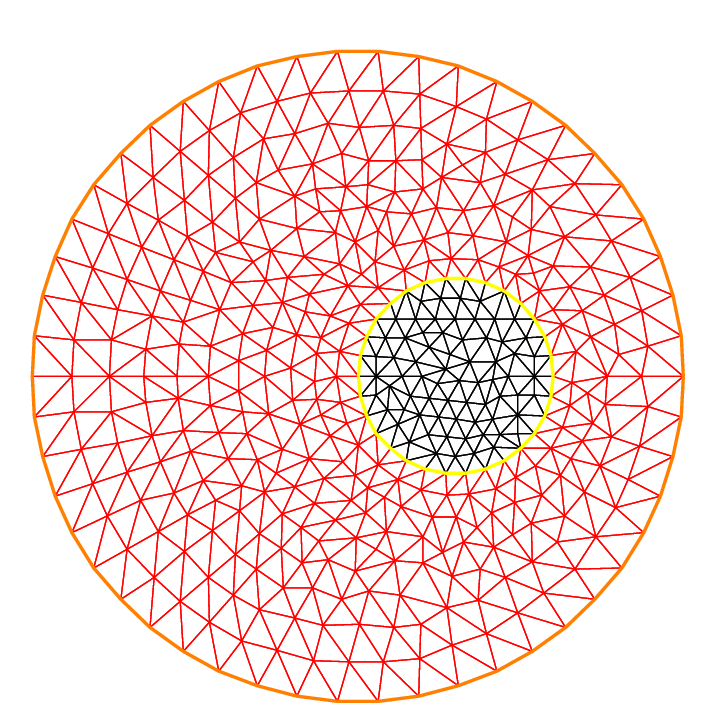"); 
plot(Thwithhole,wait=1,ps="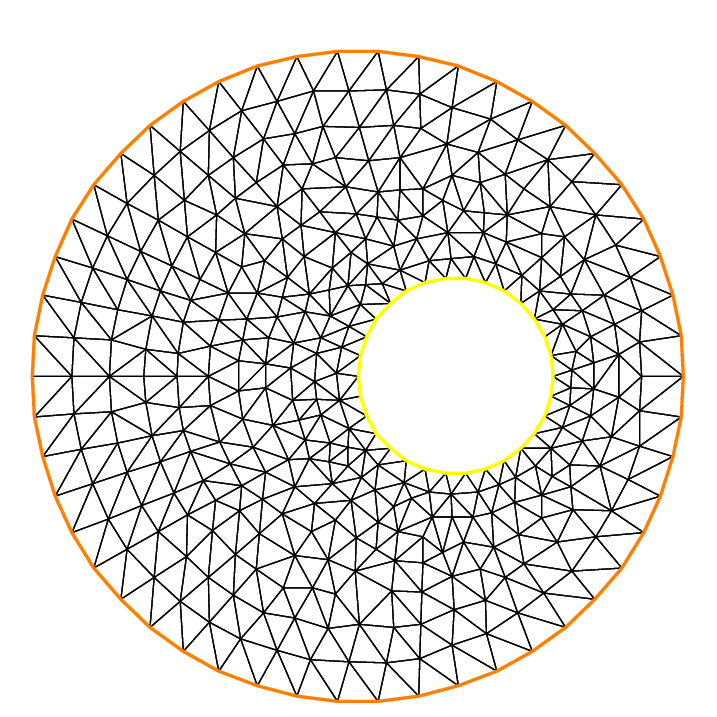"); 
\end{lstlisting}

\twoplot[height=6cm]{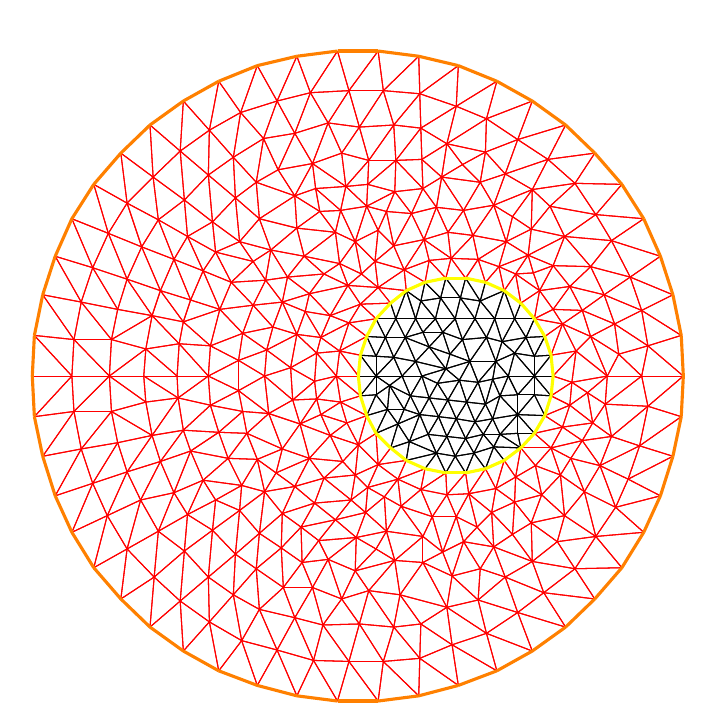}{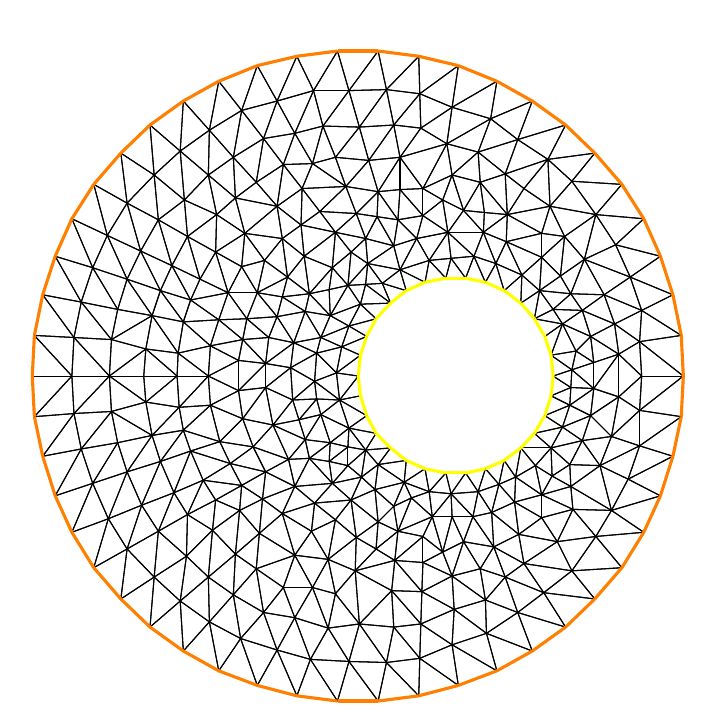}{mesh without hole}{mesh with hole}
\section{Finite Element Space}
A finite element space (F.E.S) is, usually, a space of polynomial functions on elements of $\mathcal{T}_h$, triangles here, with certain matching properties at edges, vertices, ... ; it's defined as : 
\begin{lstlisting}[firstnumber=last]
fespace Vh( Th, P1 );
\end{lstlisting}
 As of today, the known types of F.E.S. are: \textbf{P0, P03d, P1, P13d, P1dc, P1b, P1b3d, P2, P23d, P2b, P2dc, P3, P3dc, P4, P4dc, Morley, P2BR, RT0, RT03d, RT0Ortho, Edge03d, P1nc, RT1, RT1Ortho, BDM1, BDM1Ortho, TDNNS1}; where for example:
 
\begin{description}
     \item[P0,P03d]  piecewise constant discontinuous finite element  (2d, 3d), the degrees of freedom are the barycenter element value.
     \index{P0|textbf}\index{fespace!P0}
    \begin{eqnarray}
    \label{eq:P0}
     P0_{h} = \left\{ v \in L^2(\Omega) \left|\; \textrm{for all }K \in \mathcal{T}_{h}\;\;\textrm{there is }\alpha_{K}\in \R :
        \;\; v_{|K} = \alpha_{K } \right.\right\}
     \end{eqnarray}
     \item[P1,P13d]  piecewise linear  continuous finite element (2d, 3d), the degrees of freedom are the vertices values.
     \index{P1|textbf}\index{fespace!P1}\index{fespace!P13d}
     \begin{eqnarray}
     &&P1_{h} = \left\{ v \in H^{1}(\Omega) \left|\; \forall K \in \mathcal{T}_{h};
        \quad v_{|K} \in P_{1} \right.\right\} \label{eq:P1}
     \end{eqnarray}
\end{description}
We can see the description of the rest of the F.E.S. in the full documentation of \freefem.

\section{Boundary Condition}
We will see in this section how it's easy to define the boundary condition (B.C.) with \freefem, for more information about these B.C., we refer to the full documentation.
\subsection{Dirichlet B.C.}
To define Dirichlet B.C. on a border $\Gamma_d\subset \mathbb{R}$ like $u|_{\Gamma_d}=f$, we can proceed as {\ttfamily\textcolor{red}{on}(gammad,u=f)}, where {\ttfamily{u}} is the unknown function in the problem. \\
The meaning is for all degree of freedom $i$ of the boundary referred by the label ``{\ttfamily gammad}'',  the diagonal term of the matrix $a_{ii}= tgv$ with the terrible giant value $tgv$ (=$10^{30}$ by default) and the right hand side $b[i] = "(\Pi_h g)[i]" \times tgv $, where the $"(\Pi_h g)g[i]"$ is the boundary node value given by  the interpolation of $g$. (We are solving here the linear system $AX=B$, where $A=\left(a_{ij}\right)_{i=1..n;j=1..m}$ and $B=\left(b_{i}\right)_{i=1..n}$ ).\\

If $u$ is a vector like $u=(u1,u2)^T$ and we have $u1|_{\Gamma_d}=f1$ and $u2|_{\Gamma_d}=f2$, we can proceed as {\ttfamily\textcolor{red}{on}(gammad,u1=f1,u2=f2)}.
\subsection{Neumann B.C.}
The Neumann B.C. on a border $\Gamma_n\subset \mathbb{R}$, like $\ds\frac{\partial u}{\partial n}|_{\Gamma_n}=g$, appear in the Weak formulation of the problem after integrating by parts, for example $\left\langle\ds\frac{\partial u}{\partial n};\Phi\right\rangle_{\Gamma_n}=\left\langle g;\Phi\right\rangle_{\Gamma_n}=\displaystyle\int_{\Gamma_n} g\cdot \Phi\,dx=$ {\ttfamily\textcolor{red}{int1d}(Th,gamman)(g*phi)}.

\subsection{Robin B.C.}
The Robin B.C. on a border $\Gamma_r\subset \mathbb{R}$; like $a u + \kappa \ds\frac{\partial u}{\partial n} = b$ on $\Gamma_r$ where $a=a(x,y)\geq 0$, $\kappa=\kappa(x,y)\geq 0$ and $b=b(x,y)$; also appear in the Weak formulation of the problem after integrating by parts, for example $-\left\langle\kappa\ds\frac{\partial u}{\partial n};\Phi\right\rangle_{\Gamma_r}=\left\langle au-b;\Phi\right\rangle_{\Gamma_r}=\displaystyle\int_{\Gamma_r} au\cdot \Phi\,dx-\displaystyle\int_{\Gamma_r} b\cdot \Phi\,dx=$ {\ttfamily\textcolor{red}{int1d}(Th,gammar)(a*u*phi)-\textcolor{red}{int1d}(Th,gammar)(b*phi)}.\\

{\bf Important}: it is not possible to write in the same integral the linear part and the bilinear part such as in {\ttfamily\textcolor{red}{int1d}(Th,gammar)(a*u*phi-b*phi)}.

\subsection{Periodic B.C.}
In the case of Bi-Periodic B.C., they are achieved in the definition of the periodic F.E.S. such as :
\begin{lstlisting}[firstnumber=last]
fespace Vh( Th, P1,periodic=[[1,x],[3,x],[2,y],[4,y]] );
\end{lstlisting}
\section{Solve the problem}
We present here different way to solve the Poisson equation :\\
Find $u:\Omega=]0,1[\times ]0,1[\longrightarrow \mathbb{R}$ such that, for a given $f \in L^2(\Omega)$:
\begin{equation} \label{laplace}
\left\{\begin{array}{rcl}
- \Delta u &=& f \,\,  \mbox{in } \Omega\\
u&=&0 \,  \mbox{ on }\, \partial \Omega \end{array}\right.
\end{equation}
Then the basic variational formulation of (\ref{laplace}) is :\\
Find $u \in H_0^1(\Omega)$, such that for all $v\in H_0^1(\Omega)$,
\bg\label{VFlaplace}
a(u,v)=l(v)
\ed
where
$$a(u,v)=\int_{\Omega }\nabla u\cdot\nabla v\,dxdy
\mbox{ and }l(v)=\int_\Omega f\cdot v\,dxdy$$
To discretize (\ref{VFlaplace}), let $\mathcal{T}_h$ denote a regular, quasi uniform triangulation of $\Omega$ with triangles of maximum size $h<1$, let $V_h=\{v_h\in C^0(\bar{\Omega}); v_h|_T \in \mathbb{P}_1(T),\forall T\in\mathcal{T}_h;v_h=0 \mbox{ on }\partial \Omega\}$ denote a finite-dimensional subspace of $H_0^1(\Omega)$ where $\mathbb{P}_1$ is the set of polynomials of $ \mathbb{R}$ of degrees $\leq 1$.\\
Thus the discretize weak formulation of (\ref{VFlaplace}) is :
\bg \label{DVFlaplacefull}
\mbox{Find }u_h \in V_h : \int_{\Omega }\nabla u_h\cdot\nabla v_h\,dxdy-\int_\Omega f\cdot v_h\,dxdy=0\qquad \forall v_h
\in V_h. \ed
\subsection{solve}
The first method to solve (\ref{DVFlaplacefull}) is to declare and solve the problem at the same time  by using the keyword {\ttfamily\textcolor{red}{solve}} such as :

\begin{lstlisting}[firstnumber=last]
solve poisson(uh,vh,init=i,solver=LU) = // Solve Poisson Equation
	int2d(Th)( Grad(uh)'*Grad(vh) ) 	// bilinear form
	-int2d(Th)(f*vh) 					// linear form
	+on(1,2,3,4,uh=0); 					// Dirichlet B.C.
\end{lstlisting}
The solver used here is Gauss' LU factorization and when {\ttfamily init} $\neq 0$ the LU decomposition is reused so it is much faster after the first iteration. Note that if the mesh changes the matrix is reconstructed too.\\
The default solver is \textbf{sparsesolver} ( it is equal to  \textbf{UMFPACK}  if not other sparce solver is defined) or is set to
    \textbf{LU} if no direct sparse solver is available.
    The storage mode of the matrix of the underlying linear system
    depends on the type of solver chosen; for \textbf{LU}  the matrix is sky-line non
    symmetric, for \textbf{Crout} the matrix is sky-line symmetric, for
    \textbf{Cholesky} the matrix is sky-line symmetric positive
    definite,  for \textbf{CG}   the matrix is sparse symmetric positive,
    and for \textbf{GMRES}, \textbf{sparsesolver} or \textbf{UMFPACK} the matrix is just  sparse.
\subsection{problem}
The second method to solve (\ref{DVFlaplacefull}) is to declare the problem by using the keyword {\ttfamily\textcolor{red}{problem}}, and then solve it later by just call his name, such as :
\begin{lstlisting}[firstnumber=last]
problem poisson(uh,vh,init=i,solver=LU)=// Definition of the problem
	int2d(Th)( Grad(uh)'*Grad(vh) ) 	// bilinear form
	-int2d(Th)(f*vh) 					// linear form
	+on(1,2,3,4,uh=0); 					// Dirichlet B.C.
Poisson;								// Solve Poisson Equation
\end{lstlisting}
Note that, this technique is used when we have a time depend problem.
\subsection{varf}
In \freefem, it is possible to define variational forms, and use them to build matrices and vectors and store them to speed-up the script.\\

The system (\ref{DVFlaplacefull}) is equivalent to :
\bg \label{DVFlaplace}
\mbox{Find }u_h \in V_h : a(u_h,v_h)=l(v_h)\qquad \forall v_h
\in V_h. \ed
Here, 
\begin{equation}\label{defu}
u_h(x,y)=\sum_{i=0}^{M-1} {u_h}_i\phi_i(x,y)
\end{equation}
where $\phi_i={v_h}_i, i=0,...,M-1$ are the basis functions of $V_h$, $M = $ {\ttfamily Vh.ndof} is  the number of degree of freedom (i.e. the dimension of the space $V_h$) and $u_h{_i}$ is the value of $u_h$ on each degree of freedom (i.e. ${u_h}_i=${\ttfamily uh[][i]}=$U$).\\

Thus, using (\ref{defu}), we can rewrite (\ref{DVFlaplace}) such as :
\bg\label{DFVlaplacesys}
\sum_{j=0}^{M-1} A_{ij}{u_h}_j - F_i=0 ,\quad i=0,\cdots,M-1;
\ed
where
$$A_{ij}=\int_{\Omega}\nabla \phi_j\nabla\phi_i\, dxdy \qquad\mbox{ and }\qquad
F_i=\int_{\Omega}f\phi_i\, dxdy$$
The matrix $A=(A_{ij})$ is called \emph{stiffness matrix}.\\
We deduce from the above notation that (\ref{DFVlaplacesys}) is equivalent to
\bg\label{lapsys}
A\cdot U= F \Longleftrightarrow U=A^{-1}\cdot F
\ed
which can be solve in \freefem as :
\begin{lstlisting}[firstnumber=last]
int m=10,n=10;
mesh Th= square (m,n ,[x,y]);
fespace Vh( Th , P1 );
Vh uh,vh;
macro Grad (u)[dx(u),dy(u)]// in 2D

func f=1;
varf a(uh ,vh) = int2d (Th)( Grad (uh) '* Grad (vh) ) // bilinear form
+on (1 ,2 ,3 ,4 , uh =0); // Dirichlet B.C.
matrix A=a(Vh ,Vh); // build the matrix
varf l( unused ,vh) = int2d (Th)(f*vh); // linear form
Vh F; F[] = l(0, Vh); // build the right hand side vector
set (A, solver = sparsesolver );
uh [] = A^ -1*F [];
plot(uh);
\end{lstlisting}
And in 3D :
\begin{lstlisting}[firstnumber=last]
load "msh3"
load "medit"
int m=10,n=10;
mesh Th2= square (m,n ,[x,y]);
mesh3 Th=buildlayers(Th2,10,zbound=[0,1]);
fespace Vh( Th , P13d );
Vh uh,vh;
macro Grad (u)[dx(u),dy(u),dz(u)]// in 2D

func f=1;
varf a(uh ,vh) = int3d (Th)( Grad (uh) '* Grad (vh) ) // bilinear form
+on (0, 1 ,2 ,3 ,4 ,5 , uh =0); // Dirichlet B.C.
matrix A=a(Vh ,Vh); // build the matrix
varf l( unused ,vh) = int3d (Th)(f*vh); // linear form
Vh F; F[] = l(0, Vh); // build the right hand side vector
set (A, solver = sparsesolver );
uh [] = A^ -1*F [];
medit("sol",Th,uh);
\end{lstlisting}
\section{Learning by examples}
\subsection{Rate of convergence for the Poisson equation}
At the beginning, we prove that the rate of convergence in space for the Poisson equation code with $P_1$ finite element is of order 2.\\
In this example, we took zero Dirichlet homogenous B.C. on the whole boundary and we have considered the following exact solution :
$$ u_{ex}=\sin(\pi x)\cdot \sin(\pi y) $$
Then, we compute the corresponding right hand side $f(x,y)$ in order to obtain the $L^2$ norm of the error between the exact solution and the numerical one (cf. Table \ref{rate_table_poisson}) $$E(u,h_n)=|u_h(h_n)-u_{ex}(h_n)|_{L^2},\forall i=1,...,\mbox{{\ttfamily nref}},h=\delta x=1/N;$$ and then the rate of convergence in space $$r(u,h_n)=\ds\frac{\log\left(E(u,h_{n-1})/E(u,h_{n})\right)}{\log\left(h_{n-1}/h_{n}\right)},\forall i=1,...,\mbox{{\ttfamily nref}}$$
We give here a method to compute the right hand side using \textbf{Maple}\footnote{\url{http://www.maplesoft.com/}} :\\

\begin{figure}[!h]
\begin{center}
\begin{tabular}{>{\centering}m{18cm}}
\includegraphics[trim=1.9cm 15.77cm 1cm 3.3cm,clip,width=18cm]{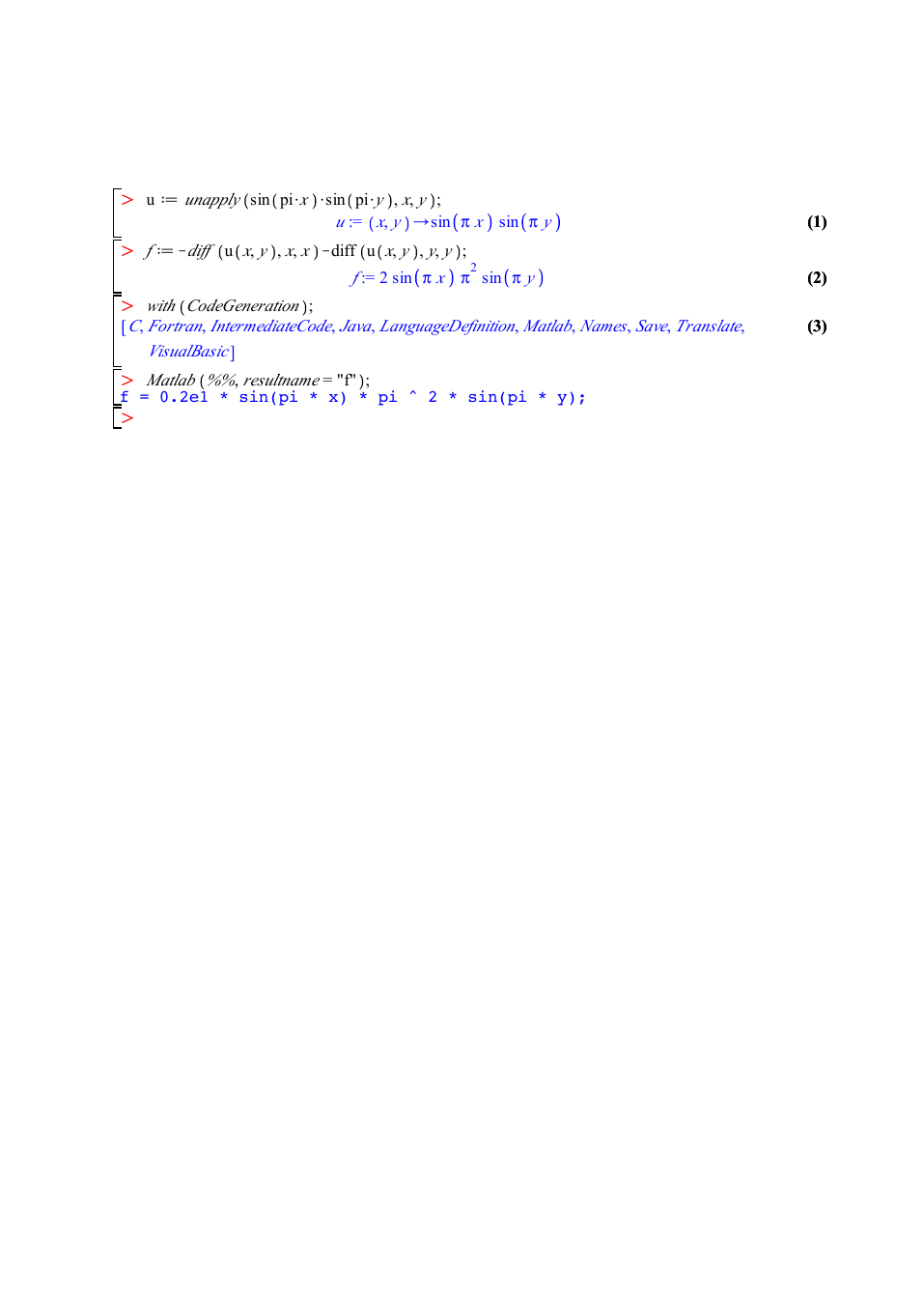}
\end{tabular}
\end{center}
\end{figure}
We can copy and paste the result of $f(x,y)$ in the \freefem code.\\
We present here the script to compute the rate of convergence in space of the code solving the Poisson equation :
\begin{lstlisting}[firstnumber=last]
int nref=4;
real[int] L2error(nref); // initialize the L2 error array
for (int n=0;n<nref;n++) {
	int N=2^(n+4); 			// space discretization
	mesh Th= square(N,N);	// mesh generation of a square
	fespace Vh(Th,P1);	// space of P1 Finite Elements
	Vh uh,vh;	// uh and vh belongs to Vh
	macro Grad(u)[dx(u),dy(u)]//
	Vh uex=sin(pi*x)*sin(pi*y); // exact solution
	Vh f=0.2e1*sin(pi*x)*pi^2*sin(pi*y); // corresponding RHS
	varf a(uh,vh) = int2d(Th)( Grad(uh)'*Grad(vh) ) // bilinear form
					+on(1,2,3,4,uh=0);	// Dirichlet B.C.
	matrix A=a(Vh,Vh);	// build the matrix
	varf l(unused ,vh) = int2d(Th)(f*vh);	//	linear	form
	Vh F; F[] = l(0,Vh);	// build the right hand side vector
	set(A,solver=sparsesolver);
	uh[] = A^-1*F[];
	L2error[n]= sqrt(int2d(Th)((uh-uex)^2));
}
for(int n=0;n<nref;n++)
	cout << "L2error " << n << " = "<<  L2error[n] <<endl;
for(int n=1;n<nref;n++)
	cout <<"convergence rate = "<< log(L2error[n-1]/L2error[n])/log(2.)  <<endl;
\end{lstlisting}
\begin{table}[!h]
\begin{center}
\begin{tabular}{|c|c|c|}
\hline
$N_n$ & $E(u,h_n)$ & $r(u,h_n)$\\
\hline
16 & 0.0047854 & - \\
\hline
32 &  0.00120952 & 1.9842\\
\hline
64 & 0.000303212 & 1.99604\\
\hline
128 & 7.58552e-05 & 1.99901\\
\hline
\end{tabular}
\end{center}
\caption{\label{rate_table_poisson} $L^2$ norm of the error and the rate of convergence.}
\end{table}%
\subsection{Poisson equation over the \textbf{Fila}'s face}
We present here a method to build a mesh from a photo using Photoshop$\circledR$ and a script in \freefem made by Frédéric Hecht.\\
We choose here to apply this method on the \textbf{Fila}'s face. To this end, we start by the photo in Figure \ref{fila_0.jpg}, and using Photoshop$\circledR$, we can remove the region that we wanted out of the domain such as in Figure \ref{fila_1}, then fill in one color your domain and use some filter in Photoshop$\circledR$ in order to smooth the boundary as in Figure \ref{fila.jpg}. Then convert your jpg photo to a pgm photo which can be read by \freefem by using in a terminal window :
\begin{lstlisting}[firstnumber=last]
convert fila.jpg fila.pgm
\end{lstlisting}
\clearpage
\threeplot[height=6cm][width=5cm,height=5.5cm][width=8cm,height=5.5cm,trim=3cm 1cm 5cm 0cm,clip]{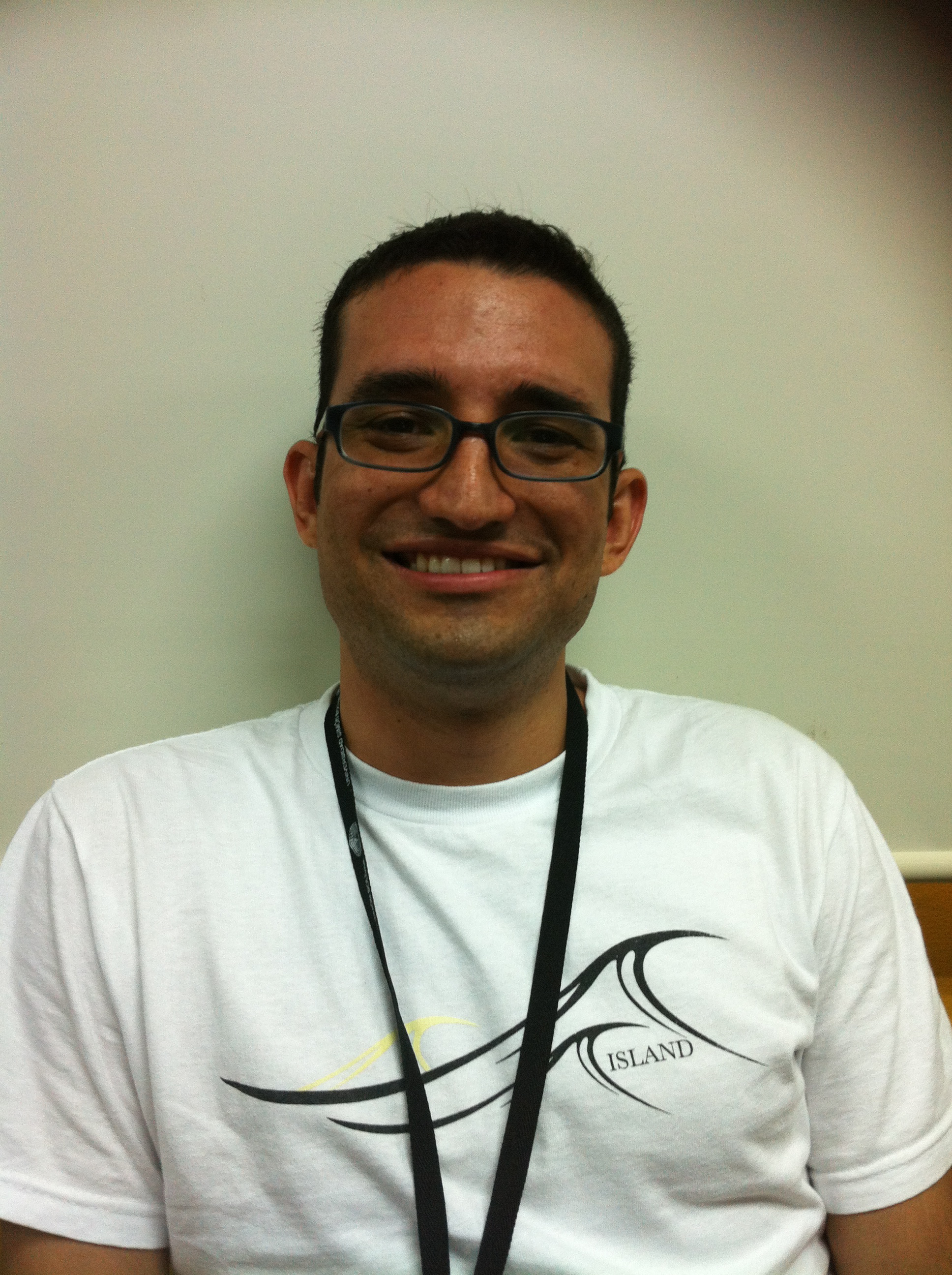}{fila_1}{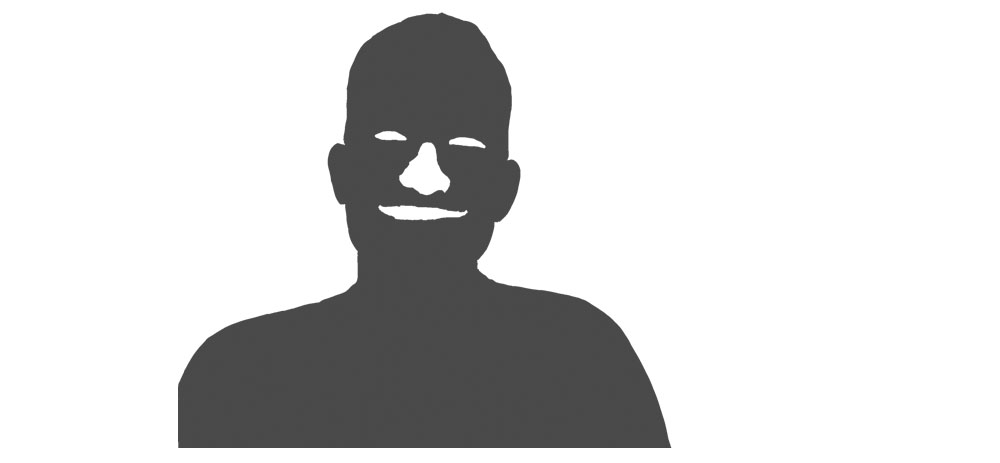}{Initial photo.}{Using Photoshop.}{Last photo.}
Finally, using the following script :
\begin{lstlisting}[firstnumber=last]
load "ppm2rnm"
load "isoline"
string fila="fila.pgm";
real[int,int] Curves(3,1);
int[int] be(1);
int nc;
{ // build the curve file xy.txt ... 
real[int,int] ff1(fila); // read  image and set to an rect. array 
//  remark (0,0) is the upper, left corner.
int nx = ff1.n, ny=ff1.m; 
// build a cartesain mesh such that the origne is qt the right place.
mesh Th=square(nx-1,ny-1,[(nx-1)*(x),(ny-1)*(1-y)]); 
// warning  the numbering is of the vertices (x,y) is 
// given by $  i = x/nx + nx* y/ny $
fespace Vh(Th,P1);
Vh f1;
f1[]=ff1; //  transforme array in finite element function.
real vmax = f1[].max ;
real vmin = f1[].min ;
real vm = (vmin+vmax)/2;
verbosity=3;
/*
 Usage of isoline
 the named parameter :
 iso=0.25   // value of iso 
 close=1, //   to force to have closing curve ... 
 beginend=be, // begin  and end of curve
 smoothing=.01, //   nb of smoothing process = size^ratio * 0.01 
 where  size is the size of the curve ...  
 ratio=0.5
 file="filename"
 
 ouptut:
  xx, yy  the array of point of the iso value 

a closed curve  number n is

in fortran notation the point of the curve are: 
(xx[i],yy[i], i = i0, i1) 
with :  i0=be[2*n],  i1=be[2*n+1];

*/
nc=isoline(Th,f1,iso=vm,close=0,Curves,beginend=be,smoothing=.005,ratio=0.1); 
verbosity=1; 
}
int ic0=be(0), ic1=be(1)-1;		
  plot([Curves(0,ic0:ic1),Curves(1,ic0:ic1)], wait=1);
// end smoothing the curve ....
macro GG(i)
border G#i(t=0,1) 
{    
  P=Curve(Curves,be(i*2),be(i*2+1)-1,t);
  label=i+1;	
} 
real lg#i=Curves(2,be(i*2+1)-1);  // 
GG(0)	GG(1)	GG(2)	GG(3)	GG(4) // number of closing curve
real hh= -10;
cout << " .. "<<endl;
func bord = G0(lg0/hh)+G1(lg1/hh)+G2(lg2/hh)+G3(lg3/hh)+G4(lg4/hh);
plot(bord,wait	=1); 
mesh Th=buildmesh(bord);
cout << " ...  "<<endl;
plot(Th,wait=1);
Th=adaptmesh(Th,5.,IsMetric=1,nbvx=1e6);
plot(Th,wait=1);
savemesh(Th,"fila.msh");
\end{lstlisting}
we can create the mesh of our domain (cf. Figure \ref{filamesh.jpg}), and then read this mesh in order to solve the Poisson equation on this domain (cf. Figure \ref{filasol.jpg})
\begin{lstlisting}[firstnumber=last]
mesh Th("fila.msh");
plot(Th);
fespace Vh(Th,P1);
Vh uh,vh;
func f = 1.;
macro Grad(u)[dx(u),dy(u)]//
solve Poisson(uh,vh) = int2d(Th)(Grad(uh)'*Grad(vh)) - int2d(Th)( f*vh) + on(1,2,3,4,5,uh=0) ;
plot(uh,dim=2,fill=true,value=true);
\end{lstlisting}
\twoplot[height=6cm]{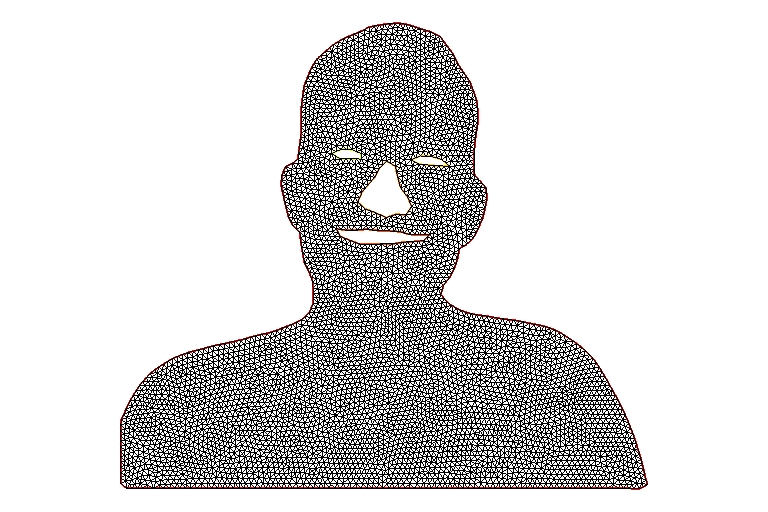}{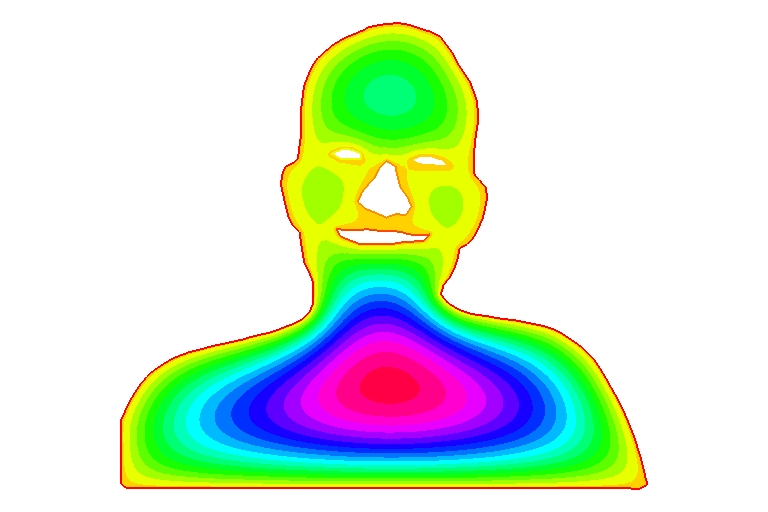}{Mesh of the \textbf{Fila}'s face.}{Solution on the \textbf{Fila}'s face.}

\subsection{Rate of convergence for an Elliptic non linear equation}
Let $\Omega=B(O,R)\subset \mathbb{R}^2$, it is proposed to solve numerically the problem which consist to find $u(x,y)$ such that
\begin{eqnarray}
\label{eqn:ellnl}
\left\{
\begin{array}{rcl}
-\Delta u(x,y) + u^3&=&f\quad \mbox{ for all }(x,y)\in\Omega\subset \R^2,
\\ 
u(x,y) &=& 0 \quad \mbox{ for all }(x,y)\mbox{ on }\p\Omega.
\end{array}\right.
\end{eqnarray}

\subsubsection{Space discretization}
Let $\mathcal{T}_h$ be the triangulation of $\Omega$ and $$V_h=\{v_h\in C^0(\bar{\Omega}); v_h|_T \in \mathbb{P}_1(T),\forall T\in\mathcal{T}_h,v_h=0 \mbox{ on }\partial \Omega\}.$$
For simplicity, we denote by $ \mathcal{V}(u)=u^2$, then the approximation of the variational formulation will be :\newline
Find $u_h\in V_h$ such that $\forall v_h\in V_h$ we have :
$$-\left\langle\Delta u_h  ;v_h\right\rangle+\left\langle  \mathcal{V}(u_h)\cdot u_h  ;v_h\right\rangle=\left\langle  f ;v_h\right\rangle$$
thus
\bg\label{FVNL}\left\langle \n u_h  ;\n v_h\right\rangle +\left\langle \mathcal{V}(u_h)\cdot u_h  ;v_h\right\rangle=\left\langle  f ;v_h\right\rangle \ed
In order to solve numerically the non linear term in (\ref{FVNL}), we will use a semi-implicit scheme such as :
$$\left\langle \n u^{n+1}_h ;\n v_h\right\rangle +\left\langle \mathcal{V}(u^n_h)\cdot u^{n+1}_h  ;v_h\right\rangle=\left\langle  f ;v_h\right\rangle,$$
and then we solve our problem by the fixed point method in this way :
$$
\begin{array}{ll}
\mbox{\tt Set } & u^0_h=u_0=0\\
\mbox{\tt Set } & \mathcal{V}(u^n_h)=\mathcal{V}(u_0)\\
\mbox{\tt Set } & err=1.\\
\mbox{\tt while $(err>1e-10)$ } & \\
&\mbox{\tt Solve } \left\langle \nabla u^{p+1}_h ;\nabla v_h\right\rangle =-\left\langle \mathcal{V}(u^n_h)\cdot u^{p}_h  ;v_h\right\rangle+\left\langle  f ;v_h\right\rangle\\
&\mbox{\tt Compute } err=\|u^n_h-u^{p+1}_h\|_{L^2}\\
&\mbox{\tt set } \mathcal{V}(u^n_h)=\mathcal{V}(u^{p+1}_h)\\
&\mbox{\tt set } u^n_h=u^{p+1}_h\\
& p=p+1;\\
 \mbox{\tt End while}&
\end{array}
$$
In order to test the convergence of this method we will study the rate of convergence in space (cf. Table \ref{rate_table_ellnl}) of the system (\ref{eqn:ellnl}) with $R=1$ and the exact solution :
$$ u_{ex}=\sin(x^2+y^2-1).$$

Then, we compute the corresponding right hand side $f(x,y)$ using \textbf{Maple} such as :\\
\begin{figure}[!h]
\begin{center}
\begin{tabular}{>{\centering}m{18cm}}
\includegraphics[trim=1.9cm 15.75cm 1cm 2.4cm,clip,width=18cm]{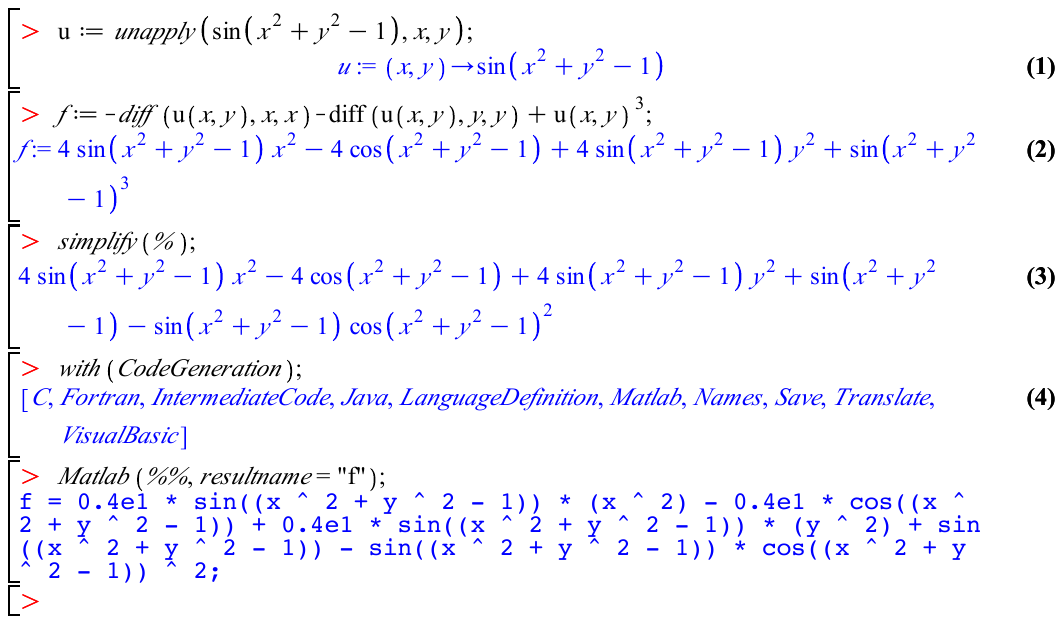}
\end{tabular}
\end{center}
\end{figure}

We present here the corresponding script to compute the rate of convergence in space of the code solving the Elliptic non linear equation (\ref{eqn:ellnl}):

\begin{lstlisting}[firstnumber=last]
verbosity=0.;
int nraff=7;
real[int] L2error(nraff); // initialize the L2 error array
for (int n=0;n<nraff;n++) {
	int N=2^(n+4); 			// space discretization
	real R=1.; // radius
	border C(t=0.,2.*pi){x=R*cos(t);y=R*sin(t);label=1;};
	mesh Th=buildmesh(C(N));
	fespace Vh(Th,P1);
	Vh uh, uh0=0, V=uh0^2, vh;
	Vh uex=sin((x ^ 2 + y ^ 2 - 1));
	Vh f=0.4e1 * sin((x ^ 2 + y ^ 2 - 1)) * (x ^ 2) - 0.4e1 * cos((x ^ 2 + y ^ 2 - 1)) + 0.4e1 * sin((x ^ 2 + y ^ 2 - 1)) * (y ^ 2) + sin ((x ^ 2 + y ^ 2 - 1)) - sin((x ^ 2 + y ^ 2 - 1)) * cos((x ^ 2 + y ^ 2 - 1)) ^ 2;
	macro Grad(u)[dx(u),dy(u)]//
	problem ELLNL(uh,vh) = 
			int2d(Th)(Grad(uh)'*Grad(vh))  // bilinear term
			+ int2d(Th) ( uh*V*vh )          // non linear term
			- int2d(Th)( f*vh )				 // right hand side
			+ on(1,uh=0);					 // Dirichlet B.C.
	real err=1.; // for the convergence
	while (err >= 1e-10){
		ELLNL;
		err=sqrt(int2d(Th)((uh-uh0)^2));
		V=uh^2;  		// actualization
		uh0=uh;
	}
	L2error[n]= sqrt(int2d(Th)((uh-uex)^2));
}
for(int n=0;n<nraff;n++)
	cout << "L2error " << n << " = "<<  L2error[n] <<endl;
for(int n=1;n<nraff;n++)
	cout <<"convergence rate = "<< log(L2error[n-1]/L2error[n])/log(2.)  <<endl;
\end{lstlisting}
\begin{table}[!h]
\begin{center}
\begin{tabular}{|c|c|c|}
\hline
$N_n$ & $E(u,h_n)$  & $r(u,h_n)$\\
\hline
16 & 0.015689 & - \\
\hline
32 &  0.0042401 & 1.88758\\
\hline
64 & 0.00117866 &   1.84695\\
\hline
128 & 0.00032964 & 1.83819\\
\hline
256 & 8.48012e-05 & 1.95873\\
\hline
512 & 1.9631e-05 & 2.11095 \\
\hline
1024 & 4.88914e-06 & 2.00548\\
\hline
\end{tabular}
\end{center}
\caption{\label{rate_table_ellnl} $L^2$ norm of the error and the rate of convergence.}
\end{table}%

\subsection{Rate of convergence for an Elliptic non linear equation with big Dirichlet B.C.}
Let $\Omega=B(O,R)\subset \mathbb{R}^2$, it is proposed to solve numerically the problem which consist to find $u(x,y)$ such that
\begin{eqnarray}
\label{eqn:ellnl}
\left\{
\begin{array}{rcl}
\Delta u(x,y) &=& \mathcal{V}(u)\cdot u\quad \mbox{ for all }(x,y)\in\Omega\subset \R^2,
\\ 
u(x,y) &=& DBC\rightarrow +\infty \quad \mbox{ for all }(x,y)\mbox{ on }\p\Omega.
\end{array}\right.
\end{eqnarray}

\subsubsection{Space discretization}
Let $\mathcal{T}_h$ be the triangulation of $\Omega$ and $$V_h=\{v_h\in C^0(\bar{\Omega}); v_h|_T \in \mathbb{P}_1(T),\forall T\in\mathcal{T}_h,v_h=p \mbox{ on }\partial \Omega\},p\longrightarrow +\infty.$$
Then the approximation of the variational formulation will be :\newline
Find $u_h\in V_h$ such that $\forall v_h\in V_h$ we have :
$$\left\langle\Delta u_h  ;v_h\right\rangle=\left\langle  \mathcal{V}(u_h)\cdot u_h  ;v_h\right\rangle$$
thus
\bg\label{FVNL}-\left\langle \n u_h  ;\n v_h\right\rangle =\left\langle \mathcal{V}(u_h)\cdot u_h  ;v_h\right\rangle \ed
In order to solve numerically the non linear term in (\ref{FVNL}), we will use a semi-implicit scheme such as :
$$-\left\langle \n u^{n+1}_h ;\n v_h\right\rangle =\left\langle \mathcal{V}(u^n_h)\cdot u^{n+1}_h  ;v_h\right\rangle,$$
and then we solve our problem by the fixed point method in this way :
$$
\begin{array}{ll}
\mbox{\tt Set } & u^0_h=u_0=DBC,p=0\\
\mbox{\tt Set } & \mathcal{V}(u^n_h)=\mathcal{V}(u_0)\\
\mbox{\tt Set } & err=1.\\
\mbox{\tt while $(err>1e-10)$ } & \\
&\mbox{\tt Solve } -\left\langle \n u^{p+1}_h ;\n v_h\right\rangle =\left\langle \mathcal{V}(u^n_h)\cdot u^{p}_h  ;v_h\right\rangle\\
&\mbox{\tt Compute } err=\|u^n_h-u^{p+1}_h\|_{L^2}\\
&\mbox{\tt set } \mathcal{V}(u^n_h)=\mathcal{V}(u^{p+1}_h)\\
& p=p+1;\\
 \mbox{\tt End while}&
\end{array}
$$
In order to test the convergence of this method we will study the rate of convergence in space (cf. Table \ref{rate_table_ellnl}) for an application of (\ref{eqn:ellnl}), where $R=1, \mathcal{V}(u)=u,DBC=0 $ or $DBC=50$.\newline
The system to be solved is then
\begin{eqnarray}
\label{eqn:ellnlrc}
\Delta u(x,y) - u^2&=& f\quad \mbox{ for all }(x,y)\in\Omega=B(O,1)\subset \R^2,
\\ \label{eqn:Dirichlet}
u(x,y) &=& DBC \quad \mbox{ for all }(x,y)\mbox{ on }\p\Omega.
\end{eqnarray}

In this case, we will the following exact solution :
$$ u_{ex}=DBC+\sin(x^2+y^2-1).$$
Then, we compute the corresponding right hand side $f(x,y)$ using \textbf{Maple} such as :\\
\begin{figure}[!h]
\begin{center}
\begin{tabular}{>{\centering}m{18cm}}
\includegraphics[trim=1.9cm 18.29cm 1cm 2.4cm,clip,width=18cm]{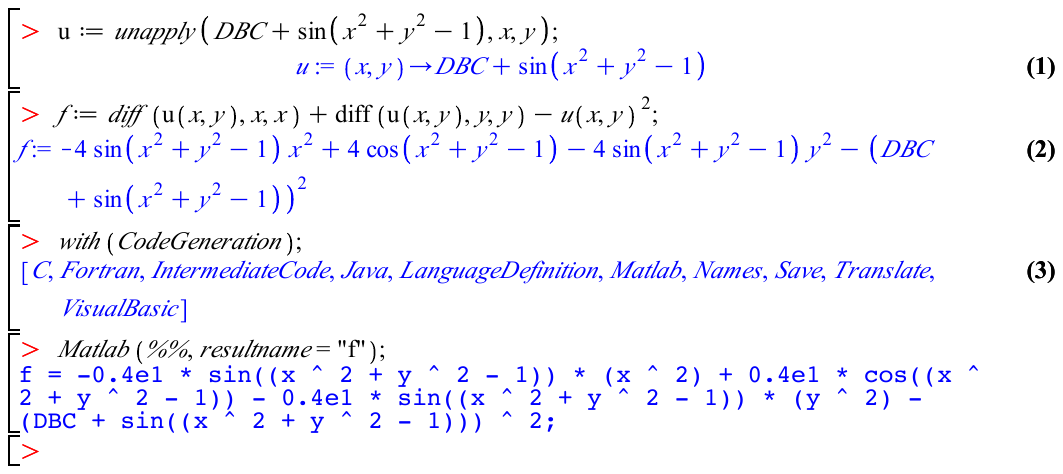}
\end{tabular}
\end{center}
\end{figure}

We present here the corresponding script to compute the rate of convergence in space of the code solving the Elliptic non linear equation (\ref{eqn:ellnlrc}):

\begin{lstlisting}[firstnumber=last]
int nref=7;
real[int] L2error(nref); // initialize the L2 error array
for (int n=0;n<nref;n++) {
	int N=2^(n+4); 			// space discretization
	real R=1.; // radius
	border C(t=0.,2.*pi){x=R*cos(t);y=R*sin(t);label=1;};
	mesh Th=buildmesh(C(N));
	fespace Vh(Th,P1);
	real DBC=0.;
	Vh uh, uh0=DBC, V=uh0, vh;
	Vh uex=DBC+sin((x ^ 2 + y ^ 2 - 1));
	Vh f=-0.4e1 * sin((x ^ 2 + y ^ 2 - 1)) * (x ^ 2) + 0.4e1 * cos((x ^ 2 + y ^ 2 - 1)) - 0.4e1 * sin((x ^ 2 + y ^ 2 - 1)) * (y ^ 2) - (DBC + sin((x ^ 2 + y ^ 2 - 1))) ^ 2;
	macro Grad(u)[dx(u),dy(u)]//
	problem ELLNL(uh,vh) = 
			- int2d(Th)(Grad(uh)'*Grad(vh))  // bilinear term
			- int2d(Th) ( uh*V*vh )          // non linear term
			- int2d(Th)( f*vh )				 // right hand side
			+ on(1,uh=DBC);					 // Dirichlet B.C.
	real err=1.; // for the convergence
	while (err >= 1e-10){
		ELLNL;
		err=sqrt(int2d(Th)((uh-V)^2));
		V=uh;  		// actualization
	}
	L2error[n]= sqrt(int2d(Th)((uh-uex)^2));
}
for(int n=0;n<nref;n++)
	cout << "L2error " << n << " = "<<  L2error[n] <<endl;
for(int n=1;n<nref;n++)
	cout <<"convergence rate = "<< log(L2error[n-1]/L2error[n])/log(2.)  <<endl;
\end{lstlisting}
\begin{table}[!h]
\begin{center}
\begin{tabular}{|c|c|c|c|c|}
\hline
$N_n$ & $E(u,h_n),$ {\ttfamily DBC=0} & $r(u,h_n)$, {\ttfamily DBC=0} & $E(u,h_n)$, {\ttfamily DBC=50} & $r(u,h_n)$, {\ttfamily DBC=50}\\
\hline
16 & 0.0159388 & - & 0.00610357 & -\\
\hline
32 &  0.00455562 & 1.80683 & 0.00244016 &1.32268\\
\hline
64 & 0.00118025 &  1.94855 & 0.000767999 &1.6678\\
\hline
128 & 0.000335335 & 1.81542 &0.000210938 &1.86429\\
\hline
256 & 8.6533e-05 & 1.95428 & 5.67798e-05 &1.89337\\
\hline
512 & 1.9715e-05 &  2.13395 &1.40771e-05 &2.01203\\
\hline
1024 & 4.90847e-06 & 2.00595 & 3.56437e-06 &1.98163\\
\hline
\end{tabular}
\end{center}
\caption{\label{rate_table_ellnl} $L^2$ norm of the error and the rate of convergence.}
\end{table}%

\subsection{Rate of convergence for the Heat equation}
Let $\Omega=]0,1[^2$, we want to solve the Heat equation :
\begin{eqnarray}
\label{eqn:heat}
\left\{
\begin{array}{rcl}
\ds\frac{\p u}{\p t}-\mu\cdot\Delta u &=& f(x,y,t)\quad \mbox{ for all }(x,y)\in\Omega\subset \R^2,t,\mu\in\mathbb{R}^+,\\
u(x,y,0) &=& u_0(x,y),\\
u&=&0 \quad\mbox{ on }\p\Omega.
\end{array}\right.
\end{eqnarray}

\subsubsection{Space discretization}
Let $\mathcal{T}_h$ be the triangulation of $\Omega$ and $$V_h=\{v_h\in C^0(\bar{\Omega}); v_h|_T \in \mathbb{P}_1(T),\forall T\in\mathcal{T}_h,v_h=0 \mbox{ on }\partial \Omega\}.$$
Then the approximation of the variational formulation will be :\newline
Find $u_h\in V_h$ such that $\forall v_h\in V_h$ we have :

$$\left\langle\frac{\p u_h}{\p t};v_h\right\rangle - \left\langle \mu\cdot\Delta u_h  ;v_h\right\rangle=\left\langle  f ;v_h\right\rangle$$
Thus
\bg\label{heatdis}
\left\langle\frac{\p u_h}{\p t};v_h\right\rangle + \mu\cdot\left\langle \n u_h  ;\n v_h\right\rangle=\left\langle  f ;v_h\right\rangle\ed

\subsubsection{Time discretization}
We will use here a $\theta$-scheme to discretize the Heat equation (\ref{heatdis}) as :
$$\left\langle \ds\frac{u^{n+1}_h-u^n_h}{\Delta t};v_h\right\rangle + \mu\cdot\theta \left\langle \n u^{n+1}_h ;\n v_h\right\rangle + \mu\cdot(1-\theta) \left\langle \n u^{n}_h ;\n v_h\right\rangle =\theta\cdot\left\langle f^{n+1};v_h\right\rangle+(1-\theta)\cdot\left\langle f^n;v_h\right\rangle.$$
Therefore :
\bg\label{heatfin}
\left\langle \ds\frac{u^{n+1}_h}{\Delta t};v_h\right\rangle + \mu\cdot\theta \left\langle \n u^{n+1}_h ;\n v_h\right\rangle = \left\langle \ds\frac{u^n_h}{\Delta t};v_h\right\rangle - \mu\cdot(1-\theta) \left\langle \n u^{n}_h ;\n v_h\right\rangle +\theta\cdot\left\langle f^{n+1};v_h\right\rangle+(1-\theta)\cdot\left\langle f^n;v_h\right\rangle. \ed
To resolve (\ref{heatfin}) with \freefem, we will write it as a linear system of the form :
$$\mathcal{A}\textbf{X}=\mathcal{B}\Longleftrightarrow\ds\sum_{j=0}^{M-1}\mathcal{A}_{i,j}\cdot\mathbf{X}_j=\mathcal{B}_{i}\mbox{ for }i=0;\dots;M-1; $$
where, $M$ is the degree of freedom, the matrix $\mathcal{A}_{i,j}$ and the arrays $\textbf{X}_{j}$ and $\mathcal{B}_{i}$ are defining as :
$$\textbf{X}_{j}=u^{n+1}_{j,h},\qquad \mathcal{A}_{i,j}=\left\{\begin{array}{ll}tgv=10^{30} &\mbox{if }i\in \p\Omega  \mbox{ and }j=i\\ \ds{\int_\Omega \frac{\varphi_i\varphi_j}{\Delta t}+\mu\cdot\theta\cdot \n \varphi_i \n \varphi_j dxdy}& \mbox{if }j\neq i
\end{array}\right.$$

$$\mathcal{B}_{i}=\left\{\begin{array}{ll}tgv=10^{30} &\mbox{if }i\in \p\Omega  \\ 
\ds{\int_\Omega \frac{ u^n_h\varphi_i}{\Delta t}-\mu\cdot(1-\theta)\cdot\n u^n_h\n\varphi_i+ \left(\theta\cdot f^{n+1}+(1-\theta)\cdot f^n\right) \varphi_i dxdy} & \mbox{otherwise}
\end{array}\right.$$

We note that the $\theta$-scheme is stable under the CFL condition (when $\theta\in[0,1/2[$) :
$$\mu\ds\frac{\Delta t}{\Delta x}+\mu\ds\frac{\Delta t}{\Delta y}\leq \ds\frac{1}{2\cdot(1-2\theta)}$$
In our test, we will consider that $\Delta x=\Delta y$ and that $CFL\in]0,1]$, then when $\theta\in[0,1/2[$, the $\theta$-scheme is stable under this condition :
$$\Delta t\leq \ds\frac{CFL\cdot(\Delta x)^2}{4\cdot \mu\cdot(1-2\theta)},$$
and for $\theta\in[1/2,1]$, the $\theta$-scheme is always stable.\\
We note also that due to the consistency error :
$$\varepsilon_{i,j}^n\leq c\Delta t\left|2\theta-1\right|+\mathcal{O}\left(\Delta x^2\right)+\mathcal{O}\left(\Delta y^2\right)+\mathcal{O}\left(\Delta t^2\right),$$
the $\theta$-scheme is consistent of order 1 in time and 2 in space when $\theta=0$ (with the CFL condition) and when $\theta=1$ (with $\Delta t=(\Delta x)^2$) and the $\theta$-scheme is consistent of order 2 in time and in space when $\theta=1/2$ (with $\Delta t=\Delta x$) (cf. Table \ref{rate_heat}).\\

\begin{rem}
When we use finite element, mass lumping is usual with explicit time-integration schemes (as when $\theta\in[0,1/2[$). It yields an easy-to-invert mass matrix at each time step, while improving the CFL condition \cite{Hug87}\label{Hug871}. In \freefem, mass lumping are defined as {\ttfamily \textcolor{red}{int2d}(Th,qft=qf1pTlump)}.  \\
\end{rem}

In order to test numerically the rate of convergence in space and in time of the $\theta$-scheme, we will consider the following exact solution :
$$ u_{ex}=\sin(\pi x)\cdot \sin(\pi y)e^{sin(t)}. $$

Then, we compute the corresponding right hand side $f(x,y)$ using \textbf{Maple} such as :\\
\clearpage
\begin{figure}[!h]
\begin{center}
\begin{tabular}{>{\centering}m{18cm}}
\includegraphics[trim=1.9cm 18.29cm 1cm 2.4cm,clip,width=18cm]{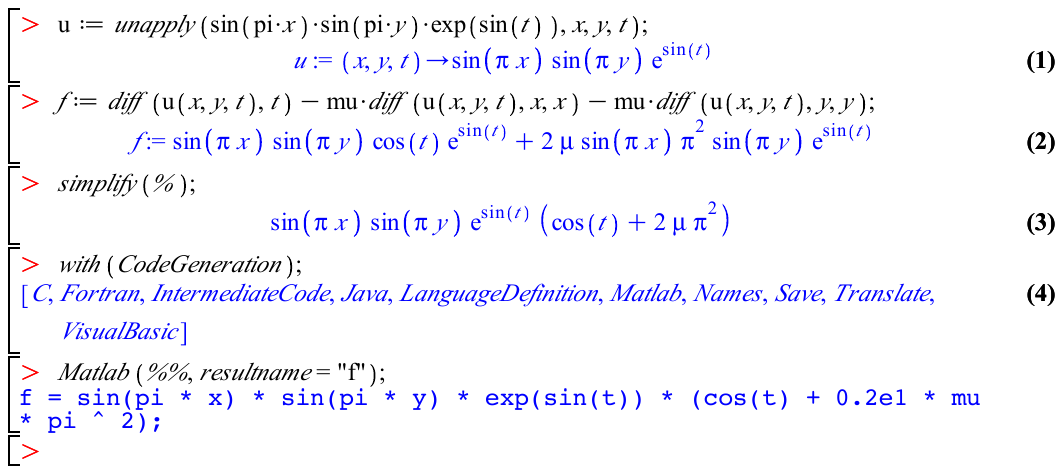}
\end{tabular}
\end{center}
\end{figure}

we remind that the rate of convergence in time is $$r(u,dt_n,h_{n})=\ds\frac{\log\left(E(u,h_{n-1})/E(u,h_{n})\right)}{\log\left(dt_{n-1}/dt_{n}\right)},\forall n=1,...,\mbox{{\ttfamily nref}}$$
\begin{lstlisting}[firstnumber=last]
macro Grad(u)[dx(u),dy(u)]//
macro uex(t) (sin(pi*x)*sin(pi*y)*exp(sin(t)))//
macro f(t)  ( sin(pi * x) * sin(pi * y) * exp(sin(t)) * (cos(t) + 0.2e1 * mu * pi ^ 2) ) //

real t, dt, h,T=.1, mu=1., CFL=1., theta=0.;
int nref=4;
real[int] L2error(nref); // initialize the L2 error array
real[int] Dx(nref); // initialize the Space discretization array
real[int] DT(nref); // initialize the Time discretization array
	
for (int n=0;n<nref;n++) {
	int N=2^(n+4);
	t=0;
	h=1./N;
	Dx[n]=h;
	if (theta<.5)
		dt=CFL*h^2/4./(1.-2.*theta)/mu; 
    else if (theta==.5)
        dt=h; 
    else
        dt=h^2;
    DT[n]=dt;
	mesh Th=square(N,N);
	fespace Vh(Th,P1);
	Vh u,u0,B;
	varf a(u,v) = int2d(Th,qft=qf1pTlump)(u*v/dt + Grad(u)'*Grad(v)*theta*mu) + on(1,2,3,4,u=0);
	matrix A = a(Vh,Vh);
	varf b(u,v) = int2d(Th,qft=qf1pTlump)(u0*v/dt - Grad(u0)'*Grad(v)*(1.-theta)*mu ) + int2d(Th,qft=qf1pTlump)( (f(t+dt)*theta+f(t)*(1.-theta))*v ) + on(1,2,3,4,u=0);
	u=uex(t);
	for (t=0;t<=T;t+=dt){
		u0=u;
		B[] = b(0,Vh);
		set(A,solver=sparsesolver);
		u[] = A^-1*B[];			
	}
	L2error[n]=sqrt(int2d(Th)(abs(u-uex(t))^2));
}
for(int n=0;n<nref;n++)
	cout << "L2error " << n << " = "<<  L2error[n] <<endl;
for(int n=1;n<nref;n++){
	cout <<"Space convergence rate = "<< log(L2error[n-1]/L2error[n])/log(Dx[n-1]/Dx[n])  <<endl;
	cout <<"Time convergence rate = "<< log(L2error[n-1]/L2error[n])/log(DT[n-1]/DT[n])  <<endl;
}
\end{lstlisting}
\begin{table}[!h]
\begin{center}
\begin{tabular}{|c|c|c|c|c|}
\hline
$N_n$ & 16 & 32 & 64 & 128\\
\hline
\hline
$E(u,h_n),\, \theta =0$& 0.00325837 & 0.000815303 & 0.000203872& 5.09709e-05\\
\hline
$r(u,h_n),\, \theta =0$ & - &  1.99874 & 1.99967 & 1.99992\\
\hline
$r(u,dt_n,h_n),\, \theta =0$ & - & 0.99937 & 0.999834 & 0.99996\\
\hline
\hline
$E(u,h_n),\, \theta =1/2$ & 0.00325537 & 0.000819141 & 0.000203817 & 5.08854e-05 \\
\hline
 $r(u,h_n),\, \theta =1/2$ & - & 1.99064 & 2.00684 & 2.00195\\
\hline
 $r(u,dt_n,h_n),\, \theta =1/2$ & - & 1.99064 & 2.00684 &  2.00195 \\
\hline
\hline
$E(u,h_n),\, \theta =1$ & 0.00323818 & 0.000807805 & 0.000201833 & 5.04512e-05 \\
\hline
 $r(u,h_n),\, \theta =1$ & - & 2.0031 & 2.00084 & 2.0002 \\
\hline
 $r(u,dt_n,h_n),\, \theta =1$ & - & 1.00155 & 1.00042 &  1.0001\\
\hline
\end{tabular}
\end{center}
\caption{\label{rate_heat} $L^2$ norm of the error and the rate of convergence in space and in time for different $\theta$.}
\end{table}%

\section{Conclusion}
We presented here a basic introduction to \freefem for the beginner with \freefem. For more information, go to the following link \url{http://www.freefem.org/ff++}. \\

{\bf Acknowledgements : }This work was done during the CIMPA School - Caracas 16-27 of April 2012. I would like to thank Frédéric Hecht (LJLL, Paris), Antoine Le Hyaric (LJLL, Paris) and Olivier Pantz (CMAP, Paris) for fruitful discussions and remarks.

\bibliographystyle{plain}

\begin{thebibliography}{}

\end{thebibliography}


\begin{thebibliography}{99}
\bibitem[Hug87]{Hug87}{\sc Thomas J. R. Hughes}. The finite element method. Prentice Hall Inc., Englewood Cliffs, NJ, 1987. Linear static and dynamic finite element analysis, With the collaboration of Robert M. Ferencz and Arthur M. Raefsky. \pageref{Hug871}\\

\bibitem[LucPir98]{LucPir98} {\sc Brigitte Lucquin and Olivier Pironneau}. {\sl Introduction to Scientific Computing}. {\em Wiley}, 1998. \href{http://www.abebooks.com/9780471972662/Introduction-Scientific-Computing-Lucquin-Brigitte-0471972665/plp}{PDF}. \pageref{LucPir981}\\
\end{thebibliography}

\end{document}